\documentclass{amsart} \textwidth=14.5cm \oddsidemargin=1cm
\usepackage{amsmath}
\usepackage{amsxtra}
\usepackage{amscd}
\usepackage{amsthm}
\usepackage{amsfonts}
\usepackage{amssymb}
\usepackage{eucal}

\newtheorem{Thm}{Theorem}
\newtheorem{Cor}{Corollary}
\newtheorem{Lem}{Lemma}
\newtheorem{Prop}{Proposition}

\theoremstyle{remark}

\newtheorem{Conj}{Conjecture}

\newcommand{\Ql}{{\overline {{\Bbb Q}_l} }}

\newcommand{\cal}{\mathcal}
\newcommand{\Gr}{{{\cal G}{\frak r} }}
\newcommand{\Fl}{{{\cal F}\ell}}
\newcommand{\tFl}{{\widetilde{{\mathcal F}\ell}}}
\newcommand{\tD}{\widetilde{D}{}}
\newcommand{\tP}{\widetilde{\mathcal P}{}}

\newcommand{\CF}{{\mathcal F}}
\newcommand{\CH}{{\mathcal H}}
\newcommand{\CS}{{\mathcal S}}
\newcommand{\CW}{{\mathcal W}}
\newcommand{\CY}{{\mathcal Y}}
\newcommand{\BW}{{\mathbb W}}
\newcommand{\BZ}{{\mathbb Z}}
\newcommand{\bm}{{\mathbf m}}
\newcommand{\fm}{{\mathfrak m}}
\newcommand{\fF}{{\mathfrak F}}
\newcommand{\fG}{{\mathfrak G}}
\newcommand{\tzeta}{{\tilde\zeta}}
\newcommand{\LGo}{{\check G}{}^0}
\newcommand{\LGsc}{{\check G}{}^{sc}}
\newcommand{\Gad}{{G^{ad}}}
\newcommand{\tFlad}{{\widetilde{{\mathcal F}\ell}{}^{ad}}}
\newcommand{\iso}{{\widetilde \longrightarrow}}

\def\square{\hbox{\vrule\vbox{\hrule\phantom{o}\hrule}\vrule}}
\newcommand{\epf}{\square}

\renewcommand{\P}{{\cal P}}

\renewcommand{\H}{{\cal H}}

\newcommand{\N}{{\cal N}}

\newcommand{\M}{{\cal M}}

\newcommand{\LG}{{\check{G}}}

\newcommand{\LT}{{\check{T}}}

\newcommand{\F}{{\cal F}}
\newcommand{\G}{{\cal G}}
\newcommand{\A}{{\cal A}}

\newcommand{\C}{{\cal C}}

\newcommand{\cons}{{\overline{\underline{{\Bbb Q}_l}}}}

\newcommand{\Z}{{\cal Z}}

\newcommand{\bI}{{\bf I}}

\newcommand{\GO}{{\bf G_O}}
\newcommand{\K}{\GO}

\newcommand{\Gm}{{\mathbb G}_m}

\newcommand{\Fq}{{\mathbb F}_q}

\newcommand{\AutO}{{\bf Aut}(O)}

\renewcommand{\proof}{{\it Proof }}

\author{Roman Bezrukavnikov, Michael Finkelberg, and Victor Ostrik}
\title
{On tensor categories attached to cells in affine Weyl groups, III}
\dedicatory{To George Lusztig with admiration}

\begin{document}
\maketitle


\begin{abstract}
We prove a weak version of Lusztig's conjecture on explicit description
of the asymptotic Hecke algebras (both finite and affine) related to
monodromic sheaves on the base affine space (both finite and affine),
and explain its relation to Lusztig's classification of character sheaves.
\end{abstract}

\begin{section}{Introduction}

\subsection{}
Truncated convolution categories were introduced by G.~Lusztig in~\cite{Lt}
as categorifications of asymptotic Hecke algebras (both in affine and
finite case). He also suggested a conjecture describing these categories
as categories of vector bundles on a square of a finite set equivariant
with respect to an algebraic group. These conjectures for affine Hecke
algebras were proved (in a weak form) in~\cite{B},~\cite{BO}.

This note is a continuation of the 
series~\cite{B},~\cite{BO}. We extend the results of {\em loc. cit.}
to the monodromic setting, and deduce the case of finite Hecke algebras.
This sheds some light on Lusztig's classification of character sheaves
on a connected reductive group in terms of representations of Drinfeld's
double of a certain finite group.

\subsection{} Let us be a bit more precise.
 For a reductive group $G$ we consider the monodromic Iwahori-equivariant
semisimple constructible complexes on the affine base affine space $\tFl$.
Here the monodromy $\zeta$ is considered as a semisimple element of the
Langlands dual group $\LG$. Their $K$-group forms an algebra
${\mathcal H}_\zeta$ with respect to convolution. This is a version of
the affine Hecke algebra. Following G.~Lusztig~\cite{Lt}, we consider certain
subcategories of semisimple monodromic perverse sheaves with truncated
convolution; their $K$-ring is a version of the asymptotic ring $J$.
Similarly to the nonmonodromic case of~\cite{BO}, we describe this
asymptotic ring in terms of convolution algebra of $F$-equivariant vector
bundles on the square of a finite set. Here $F$ is a subgroup of $\LG_\zeta$, 
the centralizer of $\zeta$ in $\LG$.

 From this description of $K$-rings on $\tFl$, we derive a similar
description of similar $K$-rings on the finite base affine space $G/U$.
Here our results are only partial, the main difficulty 
being posed by the groups
$G$ with disconnected center (see Conjecture~\ref{finite} and
Theorems~\ref{finit},~\ref{lcon}) but new even in the nonmonodromic case. 
In particular,
we prove Lusztig's conjecture~\cite{llead}~3.15 (we understand that
it was proved by Lusztig himself long ago, unpublished). We identify
a truncated convolution category on $G/U$ with ${\mathcal G}$-equivariant
vector bundles on the square of a finite set. Here ${\mathcal G}$ is the 
Lusztig's finite quotient of the algebraic group $F$. 

Lusztig has classified the isomorphism classes of
character sheaves on $G$ with the central character $\zeta$, in the
corresponding cell, in terms of his nonabelian Fourier transform
attached to ${\mathcal G}$. We relate this to the above
description of the truncated convolution category by constructing
a functor from the cell of character sheaves to the center of the
truncated convolution category.

\subsection{Acknowledgments}
R.B. was partially supported by NSF grant DMS-0505466
and  Sloan foundation award 0965-300-L216.
M.F. is grateful to V.~Serganova, I.~Zakharevich, the Weizmann Institute,
and MSRI for the hospitality and support.
He was partially supported by the CRDF award RUM1-2694,
and the ANR program ``GIMP", contract number ANR-05-BLAN-0029-01.3.
V.O. was partially supported by NSF grants DMS-0098830 and DMS-0602263.

\end{section}

\begin{section}{Monodromic central sheaves}

\subsection{}We make a free use
of the notations and results of~\cite{B},~\cite{BO}, as well as of~\cite{ab}.

We denote by $T$ the abstract Cartan group of $G$. 
We also fix a Cartan subgroup $T\subset G$. The {\em affine base
affine space} $\tFl$ is a $T$-torsor over $\Fl$. We have 
$\tFl(k)=G(F)/I_u$ where $I_u\subset I$ is the prounipotent radical.
We fix a 1-dimensional tame $\Ql$-local system $\zeta$ on $T$. 
We have $\zeta^{\otimes n}\simeq\cons$ for some integer $n$ prime to $p$.
We fix a topological generator of the tame geometric fundamental group
of $\Gm$. Then we can view
$\zeta$ as an element of $\LT(\Ql)\subset\LG(\Ql)$.
We consider the equivariant constructible derived category $\tD_I$
with respect to the adjoint action of $\bI$ on $\tFl$. We also consider
the full subcategory $\tD_I^\zeta$ of $\zeta$-monodromic complexes (with
respect to the right $T$-action on $\tFl$). Let $\tP_I\subset\tD_I,\
\tP_I^\zeta\subset\tD_I^\zeta$ be the abelian subcategories of perverse
sheaves.

The category $\tD_I^\zeta$ is monoidal; the 
convolution of $A,B\in\tD_I^\zeta$ is denoted by $A*B$.
Recall that the tensor category $Rep(\LG)$ is identified with the 
tensor category $\P_\K(\Gr)$ by means of the geometric Satake 
equivalence $S$. The construction of ~\cite{KGB} gives rise to the
{\em central} functor (see ~\cite{B}, \S 2.1) 
$\Z^\zeta:\ Rep(\LG)\to\tP_I^\zeta$.

\subsection{} 
\label{denis}
The construction of $\Z^\zeta$ proceeds as follows. We choose a smooth
curve $X$, and an $\Fq$-point $x\in X$. Following ~\cite{KGB}, \S 2.2,
we consider an indscheme $\tFl_X$ over $X$ whose $\CS$-points is the set
of quadruples $(y,\F_G,\beta,\varepsilon)$ where $y$ is an $\CS$-point of $X$;
$\F_G$ is a $G$-bundle on $X\times\CS$; furthermore, $\beta$ is a
trivialization $\F_G|_{X\times\CS-\Gamma_y}\to\F^0_G|_{X\times\CS-\Gamma_y}$
(here $\Gamma_y\subset X\times\CS$ stands for the graph of $y:\ \CS\to X$);
and finally, $\varepsilon$ is a datum of reduction of $\F_G|_{x\times\CS}$
to $U$ (here $U$ is the unipotent radical of the Borel subgroup $B\subset G$).

The argument of Proposition 3 of ~\cite{KGB} extends {\em verbatim} to 
the present situation, and proves that we have canonical isomorphisms
$\tFl_{X-x}\simeq\Gr_{X-x}\times G/U$ and $\tFl_x\simeq\tFl$.

We have a canonical embedding $T=B/U\to G/U$ (the fiber over $B/B\in G/B$).
Thus we may view the local system $\zeta$ as a perverse sheaf on 
$T\subset G/U$. Given a representation $V$ of $\LG$, the geometric Satake
equivalence produces a $\K$-equivariant perverse sheaf $S(V)\in\P_\K(\Gr)$.
By ~\cite{KGB}, \S 2.1.3, $S(V)$ gives rise to a perverse sheaf 
$S(V)_{X-x}$ on $\Gr_{X-x}$. Thus we obtain a perverse sheaf
$S(V)_{X-x}\boxtimes\zeta$ on $\tFl_{X-x}\simeq\Gr_{X-x}\times G/U$.
Taking the nearby cycles at $x\in X$ we define
$$\Z^\zeta(V):=\Psi_{\tFl_X}(S(V)_{X-x}\boxtimes\zeta)\in\P_I^\zeta.$$
The argument of ~\cite{g} extends {\em verbatim} to the present situtation
and proves that 
$\Z^\zeta:\ Rep(\LG)\to\tD_I^\zeta$ is a {\em central} functor
(see ~\cite{B}, \S 2.1) with respect to the convolution monoidal structure
on $\tD_I^\zeta$.

\subsection{}
Let $Z(\tD_I^\zeta)$ stand for the {\em center} of the monoidal category
$\tD_I^\zeta$ (see e.g. ~\cite{k}, XIII.4). Then $\Z^\zeta$ is a tensor
functor $Rep(\LG)\to Z(\tD_I^\zeta)$.

Recall that $\zeta$ can be viewed as an element of the dual torus
$\LT(\Ql)\subset\LG(\Ql)$. Let $\LG_\zeta$ be the centralizer of $\zeta$
in $\LG(\Ql)$. We have the restriction tensor functor
$Res^{\LG}_{\LG_\zeta}:\ Rep(\LG)\to Rep(\LG_\zeta)$.

\begin{Thm}
\label{main} 
There is a unique tensor functor
$\Z_\zeta:\ Rep(\LG_\zeta)\to Z(\tD_I^\zeta)$ such that
$\Z^\zeta\simeq\Z_\zeta\circ Res^{\LG}_{\LG_\zeta}:\
Rep(\LG)\to Z(\tD_I^\zeta)$.
\end{Thm}

The proof occupies subsections~\ref{fib},~\ref{mon}.

\subsection{A fiber functor}
\label{fib}
Recall that $\Lambda$ is the coweight lattice of $G$.
The argument of ~\cite{ab}, \S 3.2, provides a {\em Wakimoto sheaf}
$J_\lambda^\zeta$ in $\tD_I^\zeta$ for each $\lambda\in\Lambda$.
We denote by $\tFl_\lambda$ the preimage of the Bruhat cell
$\Fl_\lambda\subset\Fl$ under the projection $\pi:\ \tFl\to\Fl$. The locally
closed embedding $\tFl_\lambda\hookrightarrow\tFl$ is denoted by $j_\lambda$.
We fix a uniformizer $t$ of the local field $F$.
The fiber of $\pi$ over a $T$-fixed point $\lambda\in\Fl$ is canonically
a torsor over $T$, which can be canonically identified with $T$ by
the choice of a point $\lambda(t)$ in the fiber,
and there is a unique $\bI$-equivariant local system
$\zeta_\lambda$ on $\tFl_\lambda$ such that 
$\zeta_\lambda|_{\pi^{-1}(\lambda)}=\zeta$. 
For a dominant coweight $\lambda\in\Lambda^+$ we have
$J_\lambda^\zeta=j_{\lambda*}(\zeta_\lambda[\ell(\lambda)+
\operatorname{rk}_G])$. 
For an antidominant coweight $\lambda\in-\Lambda^+$ we have
$J_\lambda^\zeta=j_{\lambda!}(\zeta_\lambda[\ell(\lambda)+
\operatorname{rk}_G])$. Moreover, the proof of Theorem 5 of ~\cite{ab}
carries over {\em verbatim} to the present case, and shows that
for any $\lambda\in\Lambda$ the Wakimoto sheaf $J_\lambda^\zeta$ is
perverse, that is lies in $\P_I^\zeta$.

The proof of Proposition 5 of {\em loc. cit.} also carries over to the
monodromic case, and shows that for any $V\in Rep(\LG)$ the central sheaf
$\Z^\zeta(V)$ has a unique {\em Wakimoto filtration} 
$\BW_\bullet(\Z^\zeta(V))$ whose subquotients are the 
direct sums of Wakimoto sheaves
$J_\lambda^\zeta,\ \lambda\in\Lambda$.

Following {\em loc. cit.}, \S 3.6.5, we consider the full subcategory
$\A^\zeta\subset\P_I^\zeta$ formed by sheaves which admit a Wakimoto
filtration (whose associated graded quotient is a direct sum of
Wakimoto sheaves $J_\lambda^\zeta,\ \lambda\in\Lambda$).
We also consider the subcategory $gr\A^\zeta\subset\A^\zeta$ whose
objects are the direct sums of Wakimoto sheaves $J_\lambda^\zeta$, and
morphisms are the direct sums of isomorphisms $J_\lambda^\zeta\to
J_\lambda^\zeta$ and zero arrows. Then $\A^\zeta,gr\A^\zeta$ are monoidal
subcategories in $\tD_I^\zeta$, and $gr\A^\zeta$ is obviously equivalent
to $Rep(\LT)$. Furthermore, we see as in {\em loc. cit.} that taking the
associated graded with respect to the Wakimoto filtration is a well defined
monoidal functor $gr:\ \A^\zeta\to gr\A^\zeta\simeq Rep(\LT)$.

Since the functor $\Z^\zeta:\ Rep(\LG)\to\P_I^\zeta$ actually lands to
the subcategory $\A^\zeta$, we can compose it with $gr$ to obtain the
functor $gr\circ\Z^\zeta:\ Rep(\LG)\to gr\A^\zeta\simeq Rep(\LT)$.
The proof of Theorem 6 of {\em loc. cit.} carries over to the present
case and shows that the functor $gr\circ\Z^\zeta$ is tensor and isomorphic
to the functor of restriction to the Cartan torus $\LT\subset\LG$.

Now in order to prove theorem~\ref{main} it suffices to construct for each
$V\in Rep(\LG)$ an automorphism $a_V\in\operatorname{Aut}\Z^\zeta(V)$
such that $gr(a_V)$ is equal to the action of $\zeta\in\LT$ on 
$gr\Z^\zeta(V)$, and $a_{V_1\otimes V_2}=a_{V_1}\otimes a_{V_2}$. 
This is the subject of the next subsection.

\subsection{The monodromy transformation of the nearby cycles}
\label{mon}
All the irreducibles in $\P_I^\zeta$ are of the form 
$j_{w!*}(\zeta_w[\ell(w)+\operatorname{rk}_G])$ where 
$j_w:\ \tFl_w=\pi^{-1}(\Fl_w)\hookrightarrow\tFl$ is the locally closed
embedding of the Bruhat cell $\Fl_w$ preimage in $\tFl$, and 
$\zeta_w$ is a unique local system on $\tFl_w$ whose restriction to
$\pi^{-1}(w)\simeq T$ coincides with $\zeta$. Here $w$ is an element of
the affine Weyl group $W$. Note that $\zeta_w$ is $\bI$-equivariant
(with respect to the adjoint action) iff
$w\zeta=\zeta$. So the irreducibles in $\P_I^\zeta$ are numbered by the
stabilizer $Stab_W(\zeta)$ (where $W$ acts on $\LT(\Ql)$ through its finite
Weyl quotient $W_f$). The corresponding irreducible perverse sheaf
will be denoted by $j_{w!*}^\zeta$ for short.

Recall that the pro-algebraic group $\AutO$ of automorpshisms of $O$
acts on $\tFl$. The multiplicative group (``loop rotations'') $\Gm$
is contained in $\AutO$, and acts on $\tFl$. Any irreducible 
$j_{w!*}^\zeta$ is monodromic with respect to this $\Gm$-action.
In particular, for $\lambda\in\Lambda\subset Stab_W(\zeta)\subset W$,
the irreducible $j_{\lambda!*}^\zeta$ is $\Gm$-monodromic with the 
monodromy $\lambda^*\zeta$: here we view $\lambda$ as a cocharacter
$\Gm\to T$. The Wakimoto sheaf $J_\lambda^\zeta$ is also $\Gm$-monodromic
with the monodromy $\lambda^*\zeta$. It follows that any perverse
sheaf $\F$ in $\A^\zeta$ is $\Gm$-monodromic; its Wakimoto filtration is
preserved by the monodromy transformation $\bm$, 
and the induced transformation $gr\bm$
of $gr\F\in gr\A^\zeta\simeq Rep(\LT)$ is given by $\zeta\in\LT(\Ql)$.

For $V\in Rep(\LG)$ the nearby cycles sheaf $\Z^\zeta(V)$ comes equipped
with the monodromy action of the Galois group
$\operatorname{Gal}(F)$. It is proved in ~\cite{ab}, \S 5.2, that this
action factors through the tame fundamental group of $\Gm$, and the
corresponding monodromy transformation (recall that we have fixed a
topological generator of the tame geometric fundamental group of $\Gm$)
$\fm$ equals $\bm^{-1}$. Let $\fm_V=\fm_V^{ss}\fm_V^{un}$ be the Jordan
decomposition of the monodromy transformation of $\Z^\zeta(V)$.
It follows that the action of $gr\fm_V=gr\fm_V^{ss}$ on the fiber functor
$gr\Z^\zeta(V)\simeq Res^\LG_\LT(V)$ coincides with the action of 
$\zeta^{-1}\in\LT\subset\LG$.

Thus $a_V:=(\fm_V^{ss})^{-1}$ is the desired automorphism of $\Z^\zeta(V)$.
This completes the proof of theorem~\ref{main}. \epf

\subsection{Tilting sheaves}
\label{til}
In this and the next subsection we will show that $\Z_\zeta$ does not
factor through the restriction to a reductive subgroup which is
strictly contained in $\LG_\zeta$. This statement is not used in the
following sections.
In this subsection we assume that $G$ is adjoint, that is $\LG$ is 
simply connected (the general case will be considered in~\ref{concom}).
It is well known that in this case $\LG_\zeta$ is
connected ({\em Borel-Siebenthal subgroup}, cf. ~\cite{bs}).
Suppose that $\Z_\zeta$ factors through the restriction to a reductive subgroup
$H\subset\LG_\zeta$.
The strategy to prove the equality $H=\LG_\zeta$
is to show that the perverse sheaves $\Z^\zeta(V)$ are not ``too
decomposable'' just as sheaves, let alone central.

More precisely, for $V\in Rep(\LG)$, the perverse sheaf $\Z^\zeta(V)$
decomposes as a direct sum of indecomposable perverse sheaves, and by the
Krull-Schmidt theorem, the total number $n(\Z^\zeta(V))$ of summands
is well-defined. We can also decompose the $\LG_\zeta$-module
$Res^\LG_{\LG_\zeta}(V)$ into a direct sum of $n(V)$ irreducible
$\LG_\zeta$-modules.

The desired equality $H=\LG_\zeta$ is a consequence of the 
following

\begin{Prop} 
\label{junk}
$n(\Z^\zeta(V))=n(V)$.
\end{Prop}

In effect, the Proposition implies that any irreducible $\LG_\zeta$-module
stays irreducible after the restriction to $H$. Now apply e.g.
Lemma 8 of ~\cite{B}.
The rest of this subsection is devoted to the proof of the Proposition.

We denote by $\tD_{I,un}^\zeta$ the constructible derived category of 
$\zeta$-monodromic (with respect to the right $T$-action) $\bI$-monodromic
with unipotent monodromy (with respect to the adjoint action) complexes.
Let $\tP_{I,un}^\zeta$ be the abelian subcategory of perverse sheaves.
Recall that the irreducible perverse sheaves in $\tP_{I,un}^\zeta$
are of the form $j_{w!*}^\zeta:=j_{w!*}(\zeta_w[\ell(w)+\operatorname{rk}_G])$
where $w$ lies in the stabilizer $Stab_W(\zeta)$ of $\zeta$ in the 
(extended) affine Weyl group $W$ of $G$. We will also denote the
corresponding standard and costandard sheaves by 
$j_{w!}^\zeta,j_{w*}^\zeta$ for
short. Recall that a perverse sheaf $\F\in\tP_{I,un}^\zeta$ is called
{\em tilting} if it possesses a filtration (resp. another filtration)
with associated graded quotient equal to a direct sum of
$j_{w!}^\zeta$ (resp. $j_{w*}^\zeta$). The indecomposable tilting sheaves in
$\tP_{I,un}^\zeta$ are also parametrized by $w\in Stab_W(\zeta)$.
They are denoted by $\Xi_w^\zeta$.

Lusztig defines in ~\cite{l}, \S2, a certain subgroup $W^\zeta\subset
Stab_W(\zeta)$. This is the dual affine Weyl group (nonextended, that is
a Coxeter group) of a certain root subsystem ${\check R}{}^\zeta$ of
the root system $\check R$ of $\LG$. In fact, ${\check R}{}^\zeta$ is
the root system of $\LG_\zeta/Z(\LG_\zeta)$ (quotient modulo the center).
We have the finite Weyl subgroup $W^\zeta_f\subset W^\zeta$
(also studied by Lusztig, in ~\cite{l107}, \S1). In fact, 
$W^\zeta_f$ is the stabilizer of $\zeta\in\LT$ in the 
finite Weyl group $W_f$. Let $w^\zeta_0$ denote the longest element of
$W^\zeta_f$.

The argument of Lemma 14 of ~\cite{ab} proves that for any Iwahori-equivariant
sheaf, in particular for $\F\in\A^\zeta$
(a sheaf with a Wakimoto filtration), the convolution
$\Xi_{w^\zeta_0}^\zeta*\F$ is a perverse sheaf in $\tP_{I,un}^\zeta$.
That is the functor $?\mapsto\Xi_{w^\zeta_0}^\zeta*?,\
\tP_I^\zeta\to\tP_{I,un}^\zeta$ is exact.
It is immediate to see that this functor
annihilates no object of $\A^\zeta$.
Moreover, the argument of Theorem 7 of {\em loc. cit.} (cf. also 
Remark 9 of {\em loc. cit.}) proves that for $V\in Rep(\LG)$ the
convolution $\Xi_{w^\zeta_0}^\zeta*\Z^\zeta(V)$ is a tilting sheaf
in $\tP_{I,un}^\zeta$. So the number of indecomposable summands
$n(\Z^\zeta(V))$ is not more than the number of indecomposable summands
$n(\Xi_{w^\zeta_0}^\zeta*\Z^\zeta(V))$. 
Let $K(\tP_{I,un}^\zeta)$ denote the Grothendieck group of 
$\tP_{I,un}^\zeta$. The classes of $\Z^\zeta(V),\ \Xi_{w^\zeta_0}^\zeta$,
and of indecomposable tiltings in $K(\tP_{I,un}^\zeta)$ being known in
principle, it will be possible to conclude that
$n(\Xi_{w^\zeta_0}^\zeta*\Z^\zeta(V))=n(V)$.

More precisely, we know already that 
$[\Z^\zeta(V)]=\sum_{\lambda\in\Lambda}[\lambda:V][J_\lambda^\zeta]
=\sum_{\lambda\in\Lambda}[\lambda:V][j_{\lambda!}^\zeta]$ where
$[\lambda:V]$ is the multiplicity of $\lambda$ in $V$.
It is also known (see e.g. the next paragraph) that 
$[\Xi_{w^\zeta_0}^\zeta]=\sum_{w\in W^\zeta_f}[j_{w!}^\zeta]$.

The classes of indecomposable tiltings $[\Xi_w^\zeta]$ are obtained by
shuffling the deep results of ~\cite{l}, \cite{kt}, \cite{kt1}, ~\cite{s}. 
Namely, ~\cite{kt} relates the multiplicities $[j_{y!}^\zeta:\Xi_w^\zeta]$
to the multiplicities of Verma modules in the tilting modules over the
affine Lie algebra $\hat{\mathfrak g}$ (affinization of ${\mathfrak g}=
\operatorname{Lie}G)$ at a negative level. The latter multiplicities
are related by ~\cite{s} to the multiplicities of Verma 
$\hat{\mathfrak g}$-modules in the projectives at a positive level.
The latter multiplicities are related by ~\cite{kt1} to the stalks of
irreducible monodromic perverse sheaves on the Kashiwara affine flag scheme
$\widehat{\Fl}$ of $\hat{\mathfrak g}$. The latter stalks are governed
by the Kazhdan-Lusztig combinatorics of Lusztig's ``Master group'' $W^\zeta$
and its cosets in $W$.

Recall that $W^\zeta$ is the dual affine Weyl group of $\LG_\zeta/Z(\LG_\zeta)$.
We can consider the Kashiwara affine flag scheme
of its Langlands dual group $G'$, and the irreducible Iwahori-equivariant
perverse sheaves with trivial monodromy on this scheme.
The combinatorics of ~\cite{kt1} implies that
for $y\in W^\zeta$, the dimension of the stalks of the
corresponding irreducible perverse sheaf (for $G'$) coincides with the
dimension of the stalks of the corresponding $\zeta$-monodromic irreducible
perverse sheaf (for $G$). The argument of the previous paragraph implies
that the tilting multiplicities in
these two cases coincide as well. In the case of trivial monodromy (for $G'$),
the tilting multiplicities (in a form suited for our purposes) were computed
in ~\cite{ab}.

The bottom line is that for an irreducible $\LG_\zeta/Z(\LG_\zeta)$-module
${\mathcal V}$ with the highest weight $\mu$ the sum 
$[\Xi_{w^\zeta_0}^\zeta]*\sum_{\lambda\in\Lambda}[\lambda:{\mathcal V}]
[j_{\lambda!}^\zeta]$ equals the class $[\Xi_{w^\zeta_0\mu}^\zeta]$.
It follows that for an arbitrary irreducible $\LG_\zeta$-module ${\mathcal V}$
with the highest weight $\mu$ we also have an equality
$$[\Xi_{w^\zeta_0\mu}^\zeta]=[\Xi_{w^\zeta_0}^\zeta]*\sum_{\lambda\in\Lambda}
[\lambda:{\mathcal V}][j_{\lambda!}^\zeta]$$

We conclude that for $V\in Rep(\LG)$ we have 
$n(V)=n(\Xi_{w^\zeta_0}^\zeta*\Z^\zeta(V))\geq 
n(\Z^\zeta(V))$, and hence $n(V)=n(\Z^\zeta(V))$. Thus, $H=\LG_\zeta$
for $G$ adjoint.

\subsection{Connected components}
\label{concom}
For arbitrary $G$ let $\Gad$ denote its adjoint quotient, with the 
Langlands dual $\LGsc\to\LG$ (the simply connected cover of $\LG$).
We choose a representative $\tzeta\in\LGsc$ of $\zeta\in\LG$.
Let $\LGo_\zeta$ denote the neutral connected component of the centralizer
$\LG_\zeta$. It is known that the natural map $\LGsc_\tzeta\to\LGo_\zeta$
is surjective, and the quotient $\LG_\zeta/\LGo_\zeta$ is canonically
isomorphic to the quotient $Stab_{W_f}(\zeta)/Stab_{W_f}(\tzeta)$.

Let $\tFlad$ denote the affine base affine space for $\Gad$.
We have a natural morphism $\varpi:\ \tFl\to\tFlad$. Recall that the
connected components of $\tFl$ are numbered by the fundamental group
$\pi_1(G)$, and the connected components in the image of $\varpi$
are those in the image of $\pi_1(G)\to\pi_1(\Gad)$. Moreover, $\varpi$
is a $Z(G)$-torsor over its image.

Let $V$ be a $\LG$-module. We can also view $V$ as a $\LGsc$-module factoring
through $\LG$. We denote by $\Z^\tzeta(V)$ the corresponding central sheaf
on $\tFlad$. It is easy to see from the definitions that 
$\Z^\tzeta(V)=\varpi_*^\tzeta(\Z^\zeta(V))$ where $\varpi_*^\tzeta$ has
the following meaning. For a $\zeta$-monodromic sheaf $\F$ on $\tFl$
its direct image $\varpi_*\F$ is a direct sum of monodromic sheaves on
$\tFlad$ with monodromies ranging through the preimage of $\zeta$ in
$\LGsc$. We define $\varpi_*^\tzeta\F$ as the maximal direct summand with 
monodromy $\tzeta$. It is easy to see from the definitions that
$\varpi_*^\tzeta$ induces (the same named) tensor functor 
$Rep(H)\to Rep(\LGsc_\tzeta)$.
The tensor functor $\varpi_*^\tzeta$ induces a homomorphism of the
corresponding groups $\LGsc_\tzeta\to H$. We conclude that 
$\LGo_\zeta\subset H\subset\LG_\zeta$.

Recall that due to existence of Wakimoto filtration, the class
of a central sheaf $\F\in\Z_\zeta(Rep(H))$ in the Grothendieck group
$K(\tD_I^\zeta)={\mathbb Z}[Stab_W(\zeta)]$ actually lies inside
$\BZ[\Lambda]\subset\BZ[Stab_W(\zeta)]$. For a sheaf $\F\in\Z_\zeta(Rep(H))$
its class $[\F]\in\BZ[\Lambda]$ coincides with the character 
(i.e. restriction to the Cartan $\LT\subset H_\zeta$) of the
corresponding $H_\zeta$-module.

Now from the central property of $\F\in\Z_\zeta(Rep(H))$ we see
that the character $[\F]$ must be $Stab_{W_f}(\zeta)$-invariant.
This implies that the group of connected components
$H_\zeta/\LGo_\zeta$ must coincide with the group 
$Stab_{W_f}(\zeta)/Stab_{W_f}(\tzeta)=\LG_\zeta/\LGo_\zeta$.
We conclude that $H=\LG_\zeta$.

\end{section}

\begin{section}{Generalized central sheaves}

\subsection{Generalized Wakimoto sheaves}
Recall that $Stab_W(\zeta)$ is isomorphic to 
the semidirect product $W(\LG_\zeta^0)\rtimes
(\LG_\zeta/\LG_\zeta^0)$. The subgroup $\LG_\zeta/\LG_\zeta^0\subset
Stab_W(\zeta)$ actually lies
inside $Stab_{W_f}(\zeta)$ in the finite Weyl group, and coincides with
the subgroup of elements of length zero in $Stab_{W_f}(\zeta)$. We will
denote this length zero subgroup of $Stab_{W_f}(\zeta)$ by $W_\zeta^0$
for brevity. Note that for $w\in W_\zeta^0$ we have
$j^\zeta_{w!}=j^\zeta_{w!*}=j^\zeta_{w*}$, and this clean sheaf is supported at the finite
base affine space $G/U\subset\tFl$.

The subgroup $W_\zeta^0\subset W_f$ acts on the lattice $\Lambda$,
and we denote by $\widehat\Lambda$ the semidirect product 
$\Lambda\rtimes W_\zeta^0$.
For $\widehat{\lambda}=(\lambda,w)\in\widehat\Lambda$ we consider the
{\em generalized Wakimoto sheaf} $J^\zeta_{\widehat\lambda}=
J^\zeta_\lambda*j^\zeta_{w!}$. We have $J^\zeta_{\widehat\lambda}*
J^\zeta_{\widehat\mu}\simeq J^\zeta_{\widehat{\lambda}\cdot\widehat\mu}$.
We consider the category $gr\widehat{\A}{}^\zeta$ formed by the direct
sums of generalized Wakimoto sheaves, whose morphisms are the direct
sums of isomorphisms $J^\zeta_{\widehat\lambda}\to J^\zeta_{\widehat\lambda}$
and zero arrows. It is a monoidal subcategory of $\tD_I^\zeta$.

\begin{Lem}
\label{Eugene}
The monoidal category $gr\widehat{\A}{}^\zeta$ is equivalent to
the category of coherent sheaves on $\widehat\Lambda$ with convolution.
\end{Lem}

\proof: Recall that $J^\zeta_{\widehat\lambda}*
J^\zeta_{\widehat\mu}\simeq J^\zeta_{\widehat{\lambda}\cdot\widehat\mu}$,
however we need a more precise choice of this isomorphism. To this end
we will choose a special subgroup 
$\widetilde\Lambda$ of the preimage of $\widehat\Lambda$ in
the normalizer $N_F(T)$ of the torus $T$ in $G(F)$. To begin with, we will
choose a special subgroup 
$\widetilde{W}{}_\zeta^0$ of the preimage of $W_\zeta^0$ 
in the normalizer of $T$
in $G$. Namely, we need the kernel $K_\zeta$ of the projection
$\widetilde{\Lambda}\twoheadrightarrow\widehat\Lambda$ (the same as the
kernel of the projection 
$\widetilde{W}{}_\zeta^0\twoheadrightarrow W_\zeta^0$) to be finite and
{\em central} in
$\widetilde{\Lambda}$ (and in $\widetilde{W}{}_\zeta^0$).

Recall that $G$ is almost simple. It implies that the center of $G$ is
cyclic except for the case $G=Spin_{4n}$. In case $G=Spin_{4n}$, we have
$\LG=PSO_{4n}$, and it is easy to see that the whole Weyl group of 
$\LG$ can be lifted (isomorphically) into the adjoint group $PSO_{4n}$. 
In particular, $W_\zeta^0\subset W_f$ is lifted (isomorphically) into
$PSO_{4n}$.
We define $\widetilde{W}{}_\zeta^0$ as the full preimage of this lift in the
simply connected cover $Spin_{4n}\twoheadrightarrow PSO_{4n}$.
Together with
the canonical lift of $\Lambda$ into $N_F(T)$, it generates the desired 
subgroup
$\widetilde{\Lambda}\twoheadrightarrow\widehat\Lambda$.

In case the center of $G$ is cyclic, the group $W_\zeta^0$ is cyclic as well,
say, generated by $w$. Let us choose an arbitrary lift $\widetilde{w}\in
N(T)$ of finite order. Then $\widetilde{w}$ generates the desired
subgroup $\widetilde{W}{}_\zeta^0$. Together with
the canonical lift of $\Lambda$ into $N_F(T)$, it generates the desired lift
$\widetilde{\Lambda}\twoheadrightarrow\widehat\Lambda$.

We can view $\widetilde{\lambda}\in\widetilde\Lambda$ as a point of $\tFl$
lying over the point $\widehat{\lambda}\in W\subset\Fl$.
Let $J^\zeta_{\widetilde\lambda}$ stand for the Wakimoto sheaf 
$J^\zeta_{\widehat\lambda}$ equipped with a trivialization of its stalk
at the point $\widetilde\lambda$. Then we have a canonical isomorphism
$J^\zeta_{\widetilde\lambda}*
J^\zeta_{\widetilde\mu}\cong J^\zeta_{\widetilde{\lambda}\cdot\widetilde\mu}$.

Let $Coh(\widetilde{\Lambda})$ be the category of coherent sheaves with finite
support on the discrete set $\widetilde{\Lambda}$. The category $Coh(\widetilde{\Lambda})$
is monoidal via the convolution induced by the multiplication in the group $\widetilde{\Lambda}$.
We have the functor $F: Coh(\widetilde{\Lambda})\to gr\widehat{\A}{}^\zeta$ defined
via $F(X)=\oplus_{\widetilde{\lambda}\in \widetilde{\Lambda}}X_{\widetilde{\lambda}}\otimes
J^\zeta_{\widetilde{\lambda}}$. The isomorphisms from the previous paragraph determine
the tensor structure on the functor $F$. Note that the functor $F$ is surjective in a sense that
any object of $gr\widehat{\A}{}^\zeta$ is a direct summand of the object of the form $F(?)$;
in this sense the category $gr\widehat{\A}{}^\zeta$ is a quotient of the category 
$Coh(\widetilde{\Lambda})$. We can describe this situation in the following way.
Let $I: gr\widehat{\A}{}^\zeta\to Coh(\widetilde{\Lambda})$ be the right adjoint functor of $F$. 
The object $A:=I({\bf 1})\in Coh(\widetilde{\Lambda})$ has a natural structure of the object
of Drinfeld center of the category $Coh(\widetilde{\Lambda})$ given by the isomorphisms
$c_A: I({\bf 1})*X\simeq I(F(X))\simeq X*I({\bf 1})$. Moreover, $A$ is naturally a commutative
algebra in the Drinfeld center of the category $Coh(\widetilde{\Lambda})$. Now the monoidal 
category $gr\widehat{\A}{}^\zeta$ can be identified with the category of $A-$modules in the category
$Coh(\widetilde{\Lambda})$ with monoidal structure given by the tensor product over $A$.
Thus to identify the monoidal category  $gr\widehat{\A}{}^\zeta$ we just need to identify
the object $A\in Coh(\widetilde{\Lambda})$ as a commutative algebra in the Drinfeld center
of $Coh(\widetilde{\Lambda})$.

Observe that the functor $I$ has the following description. Any object $X$ of the category
$gr\widehat{\A}{}^\zeta$ admits a canonical decomposition $X=\oplus_{\widehat\lambda \in
\widehat\Lambda}X_{\widehat\lambda}$ where each $X_{\widehat\lambda}$ is a direct sum of several
copies of $J^\zeta_{\widehat\lambda}$. Now $I(X_{\widehat\lambda})$ is the sheaf obtained by
restriction of $X_{\widehat\lambda}$ to the fiber over $\widehat\lambda$ of the projection
$\widetilde{\Lambda}\to \widehat\Lambda$ (we identify here coherent and constructible sheaves
over a finite set). This together with the fact that the kernel of $\widetilde\Lambda \to \widehat\Lambda$
is central implies that the following diagram commutes:

\begin{equation}
\begin{CD} 
A*X@>c_A>>X*A\\
@VV\Gamma V @VV\Gamma V\\
\Gamma(A)\otimes \Gamma(X)@>c>>\Gamma(X)\otimes \Gamma(A)
\end{CD}
\end{equation}
where $\Gamma : Coh(\widetilde{\Lambda})\to {\mbox Vec}$ is the tensor functor of global sections
and $c$ is the usual commutativity morphism in the category of vector spaces ${\mbox Vec}$.
It is easy to see that this condition identifies $A$ uniquely as an object of the Drinfeld center
of the category $Coh(\widetilde{\Lambda})$. Moreover it is obvious that the object $A$ has
unique structure of semisimple commutative algebra in this category.

Now we note that the description of the object $A$ above coincides with the description of a
similar object $A'$ associated with obvious tensor functor $Coh(\widetilde{\Lambda})\to
Coh(\widehat{\Lambda})$. The Lemma is proved.

\epf

\subsection{Remark} 
Recall that the category $Coh({\mathbb Z}/2{\mathbb Z})$ admits 
two different associativity constraints: the usual one and the twisted one which
differs from the first by an element from $H^3({\mathbb Z}/2{\mathbb Z}, \Gm)$
(see e.g. \cite{ENO}). The difficulty in the proof of Lemma 1 is caused by the
fact that both these categories are quotients of the category $Coh({\mathbb Z}/4{\mathbb Z})$ 
(with the usual associativity constraint). This is a consequence of the fact that the
map $H^3({\mathbb Z}/2{\mathbb Z}, \Gm)\to H^3({\mathbb Z}/4{\mathbb Z}, \Gm)$
induced by the projection ${\mathbb Z}/4{\mathbb Z}\to {\mathbb Z}/2{\mathbb Z}$
is zero.

\subsection{Decentralized central sheaves}
Recall the central functor $\Z^\zeta:\ Rep(\LG)\to Z(\tD^\zeta_I)$.
Composing it with the forgetting of the central structure we get the
monoidal functor $\Z^\zeta_0:\ Rep(\LG)\to\tD^\zeta_I$. 
The proof of Theorem~\ref{main} admits the following corollary:

\begin{Cor}
\label{msri}
There is a unique monoidal functor
$\Z_\zeta^0:\ Rep(\LG_\zeta^0)\to\tD^\zeta_I$ such that
$\Z^\zeta_0\simeq\Z_\zeta^0\circ Res^{\LG}_{\LG_\zeta^0}:\
Rep(\LG)\to\tD_I^\zeta$. 
\end{Cor}

\proof: In case $G$ is adjoint there is nothing to prove since
$\LG_\zeta=\LG_\zeta^0$. The general case is deduced from the 
adjoint case as follows.
Recall the setup and notations of~\ref{concom}.
We have the following diagram of monoidal functors
\begin{equation}
\begin{CD} 
Rep(\LG_\zeta)@>\Z_\zeta>>\tD_I^\zeta\\
@VV\operatorname{res}V @VV\varpi_*^{\tilde\zeta}V\\
Rep(\LGsc_{\tilde\zeta})@>\Z_{\tilde\zeta}>>\tD_I^{\tilde\zeta}
\end{CD}
\end{equation}

The restriction $\operatorname{res}$ factors through
$Rep(\LG_\zeta)\to Rep(\LG_\zeta^0)\to Rep(\LGsc_{\tilde\zeta})$.
Now note that $Rep(\LGsc_\tzeta)$ is decomposed into blocks according
to the characters of the central subgroup $Z(\LGsc)\cap\LGsc_\tzeta$,
and $\Z_\tzeta$ is compatible with the decomposition of
$\tD_I^\tzeta$ into blocks corresponding 
to the connected components of $\tFlad$. Now the block corresponding
to the trivial character is $Rep(\LG_\zeta^0)$, while 
$\varpi_*^\tzeta$ sends $\tD_I^\zeta$ into the blocks corresponding
to the connected components of $\tFl$.
The Corollary follows. \epf

\subsection{Generalized central sheaves}
We consider the semisimple abelian category $\CW_\zeta^0$ formed by the
direct sums of perverse sheaves $j^\zeta_{w!*},\ w\in W_\zeta^0$, on $G/U$.
It is a monoidal subcategory of $\tD_I^\zeta$ with respect to convolution.
According to Lemma~\ref{Eugene}, $\CW_\zeta^0$ is monoidally equivalent
to $Coh(W_\zeta^0)$. We will denote by 
$\varsigma:\ Coh(W_\zeta^0)\iso\CW_\zeta^0$
the monoidal equivalence of Lemma~\ref{Eugene}.
We consider the external tensor product of monoidal categories
$Rep(\LG)\boxtimes Coh(W_\zeta^0)$. It is a monoidal category equivalent
to the category $Coh(\LG\backslash(\LG\times W_\zeta^0)/\LG)$ of 
left and right $\LG$-equivariant coherent sheaves on $\LG\times W_\zeta^0$
(with monoidal structure given by convolution).

The functor $\Z^\zeta:\ Rep(\LG)\to\P^\zeta_I$ extends to the monoidal
functor $\widehat\Z{}^\zeta:\ Coh(\LG\backslash(\LG\times W_\zeta^0)/\LG)=
Rep(\LG)\boxtimes Coh(W_\zeta^0)\to\tD^\zeta_I$ defined as follows
(notations of~\ref{denis}):
$$\widehat\Z{}^\zeta(V\boxtimes\CF):=\Psi_{\tFl_X}(S(V)_{X-x}\boxtimes 
\varsigma(\CF))\in\P_I^\zeta.$$
Now recall that the group of connected components of $\LG_\zeta$ is
canonically isomorphic to $W_\zeta^0$. Thus we have the diagonal embedding
$\Delta:\ \LG_\zeta\hookrightarrow\LG\times W_\zeta^0$. Evidently,
$\Delta^{-1}(\LG\times\{e\})=\LG_\zeta^0$. This embedding gives rise to
the monoidal functor
$Res_\Delta:\ Coh(\LG\backslash(\LG\times W_\zeta^0)/\LG)\to
Coh(\LG_\zeta^0\backslash\LG_\zeta/\LG_\zeta^0)$.

\begin{Prop}
\label{key}
There is a unique monoidal functor $\widehat{\Z}{}^0_\zeta:\
Coh(\LG_\zeta^0\backslash\LG_\zeta/\LG_\zeta^0)\to\tD_I^\zeta$
such that $\widehat\Z{}^\zeta\simeq\widehat\Z{}^0_\zeta\circ Res_\Delta:\
Coh(\LG\backslash(\LG\times W_\zeta^0)/\LG)\to\tD^\zeta_I$.
\end{Prop}

\proof: We have a full monoidal subcategory 
$Rep(\LG^0_\zeta)=Coh(\LG^0_\zeta\backslash\LG^0_\zeta/\LG^0_\zeta)\subset
Coh(\LG^0_\zeta\backslash\LG_\zeta/\LG_\zeta^0)$. An object $X$ of
$Coh(\LG^0_\zeta\backslash\LG_\zeta/\LG_\zeta^0)$ decomposes as a
direct sum $\bigoplus_i U_i*Res_\Delta(V_i)$ for certain 
$U_i\in Rep(\LG^0_\zeta)$, and $V_i\in Rep(\LG)\boxtimes Coh(W_\zeta^0)=
Coh(\LG\backslash(\LG\times W_\zeta^0)/\LG)$.
We define $\widehat{\Z}{}^0_\zeta(X):=
\bigoplus_i\Z^0_\zeta(U_i)*\widehat\Z{}^\zeta(V_i)$
(for the notation $\Z^0_\zeta$ see Corollary~\ref{msri}).
It is immediate to check that $\widehat\Z{}^0_\zeta$ is well defined
and has a monoidal structure.
\epf

\subsection{Remark}
\label{some}
The group $\LG_\zeta$ acts naturally on $W_\zeta^0$. It also acts on 
$W_\zeta^0\times W_\zeta^0$ diagonally, and we can consider the category
$Coh_{\LG_\zeta}(W_\zeta^0\times W_\zeta^0)$ of $\LG_\zeta$-equivariant
coherent sheaves on the finite set $W_\zeta^0\times W_\zeta^0$. This is a
monoidal category with respect to convolution. Clearly, 
$Coh_{\LG_\zeta}(W_\zeta^0\times W_\zeta^0)$ is monoidally equivalent to 
$Coh(\LG_\zeta^0\backslash\LG_\zeta/\LG_\zeta^0)$. So the above Proposition
provides us with the monoidal functor $\widehat{\Z}{}^0_\zeta:\
Coh_{\LG_\zeta}(W_\zeta^0\times W_\zeta^0)\to\tD_I^\zeta$.

\end{section}

\begin{section}{Cells in affine Weyl groups}

\subsection{Extended affine Hecke algebras}
We are going to define the {\em extended dual affine Hecke algebra}
${\mathcal H}_\zeta={\mathcal H}(\LG_\zeta)$ of the 
(not necessarily connected) reductive group
$\LG_\zeta$. Let $\CH_\zeta^0$ stand for the dual affine Hecke algebra
$\CH(\LG_\zeta^0)$
of the connected reductive group $\LG_\zeta^0$ (same as affine Hecke algebra
of the reductive group $'G$ Langlands dual to $\LG_\zeta^0$)\footnote{It is
defined in a standard way in terms of the dual affine Weyl group
$W(\LG_\zeta^0)$
of $\LG_\zeta^0$ (same as the affine Weyl group of $'G$) and its length 
function; which are defined in~\cite{cells4},~1.6 
for arbitrary connected reductive groups.}; it is a 
${\mathbb Z}[v,v^{-1}]$-algebra with the Kazhdan-Lusztig basis
$\{C_y,\ y\in W(\LG_\zeta^0)\}$. The group of connected components
$W_\zeta^0=\LG_\zeta/\LG_\zeta^0$ acts by outer automorphisms of $\LG_\zeta^0$,
and hence it acts by the automorphisms of $W(\LG_\zeta^0)$ preserving
the length function. We denote by $\Upsilon:\ W_\zeta^0
\to Aut(W(\LG_\zeta^0))$ the corresponding homomorphism, and keep the
same notation for the corresponding homomorphism 
$\Upsilon:\ W_\zeta^0\to Aut(\CH(\LG_\zeta^0))$ preserving the
Kazhdan-Lusztig basis.


Note that $Stab_W(\zeta)$ is the semidirect product $W(\LG_\zeta^0)\rtimes
W_\zeta^0$.
We define ${\mathcal H}_\zeta$ as the
semidirect product of $\CH(\LG_\zeta^0)$ with $W_\zeta^0$,
that is tensor product 
$\CH(\LG_\zeta^0)\otimes_{\mathbb Z}{\mathbb Z}[W_\zeta^0]$
with multiplication $(h_1\otimes x_1)\cdot(h_2\otimes x_2)=
h_1\cdot\Upsilon(x_1)h_2\otimes x_1\cdot x_2$. For $w=(y,x)\in 
Stab_W(\zeta)=W(\LG^0_\zeta)\rtimes W_\zeta^0$ 
we define the {\em Kazhdan-Lusztig element} $C_w\in{\mathcal H}_\zeta$ as
$C_y\otimes x$ where 
$C_y\in\CH(\LG_\zeta^0)$ is an element of the Kazhdan-Lusztig basis.

\subsection{} We consider the category ${\mathfrak C}_\zeta$ whose
objects are the finite direct sums of simple perverse sheaves
$j_{w!*}^\zeta$ with shifts (here $w\in Stab_W(\zeta)$); the morphisms
are those in the derived category. It is a monoidal category with
respect to the convolution. Let ${\mathcal K}_\zeta$ be the abelian group with
generators corresponding to the isomorphism classes of objects of
${\mathfrak C}_\zeta$, and relations $[z_1]=[z_2]+[z_3]$ whenever $z_1$ is
isomorphic to $z_2\oplus z_3$. We regard ${\mathcal K}_\zeta$ as a
${\mathbb Z}[v,v^{-1}]$-module by $v^mc=c[-m]$. The convolution on 
${\mathfrak C}_\zeta$ gives rise to the structure of 
${\mathbb Z}[v,v^{-1}]$-algebra on ${\mathcal K}_\zeta$. We consider a
homomorphism $r$ of ${\mathbb Z}[v,v^{-1}]$-modules from 
${\mathcal H}_\zeta$ to ${\mathcal K}_\zeta$ sending $C_w$ to
$[j_{w!*}^\zeta]$ for $w\in Stab_W(\zeta)$. The following proposition
is essentially a reformulation of~\cite{l}, Proposition~5.4. The proof is
parallel to the case of $G/B$ considered in~\cite{l1},~13.2.

\begin{Prop}
\label{skip}
$r$ is an isomorphism of ${\mathbb Z}[v,v^{-1}]$-algebras.
\end{Prop}
\epf

\subsection{Truncated convolution category} 
\label{J}
Various notions and facts
about the affine Hecke algebras can be carried over to the
algebra ${\mathcal H}_\zeta$. In particular the concept of
{\em cells} (left, right and two-sided) is defined in ${\mathcal H}_\zeta$
(with respect to the basis $C_w$), see~\cite{cells4},~1.1. 
Similarly, the $a$-function $a: Stab_W(\zeta)\to {\mathbb Z}$,
and the distinguished involutions are defined, and each left (or right) 
cell contains a unique distinguished involution, see {\em loc. cit.}
It is known that the $a$-function is constant on the two-sided cells
and hence the expression $a({\underline c})$ makes sense. The group
$W^0_\zeta$ acts on the set of two-sided cells for the
Hecke algebra ${\mathcal H}(\LG_\zeta^0)$. It is easy to see that
the two-sided cells for ${\mathcal H}_\zeta$ are in the bijection with the
$W^0_\zeta$-orbits on the cells of ${\mathcal H}(\LG_\zeta^0)$.
Combining this with the Lusztig's bijection between the two-sided cells
for ${\mathcal H}(\LG_\zeta^0)$ and the unipotent classes for $\LG_\zeta^0$ 
we obtain that 
the set of two-sided cells for ${\mathcal H}_\zeta$ is in natural bijection
${\underline c}\mapsto u({\underline c})$ with the set of unipotent 
classes in the group $\LG_\zeta$. 
Finally, similarly to the case of usual affine Hecke algebra for any
two-sided cell ${\underline c}$ one defines the canonical right cell
$\Gamma_{\underline c} \subset {\underline c}$, see~\cite{LX}.

For a two-sided cell ${\underline c}$ let $Stab_W(\zeta)_{\le {\underline c}}=
\bigcup_{{\underline c}'\le_{LR} {\underline c}}{\underline c}'$ and
$Stab_W(\zeta)_{< {\underline c}}=\bigcup_{{\underline c}'<_{LR} 
{\underline c}}{\underline c}'$. 
Let $\C_\zeta$ be the semisimple abelian subcategory of $\tD^\zeta_I$
formed by the direct sums of irreducible perverse sheaves
$j_{w!*}^\zeta,\ w\in Stab_W(\zeta)$.
Let ${\mathcal C}_\zeta^{\le {\underline c}}$
(respectively ${\mathcal C}_\zeta^{< {\underline c}}$) be the Serre 
subcategory of ${\mathcal C}_\zeta$ generated by the objects $j_{w!*}^\zeta$
where $w\in Stab_W(\zeta)_{\le {\underline c}}$ (respectively 
$w\in Stab_W(\zeta)_{< {\underline c}}$). Consider the Serre quotient
category ${\mathcal C}_\zeta^{\underline c}= {\mathcal C}_\zeta^{\le 
{\underline c}}/{\mathcal C}_\zeta^{< {\underline c}}$. The category
${\mathcal C}_\zeta^{\underline c}$ is endowed with the structure of
monoidal category via {\em truncated convolution} $X\bullet Y=
{}^pH^{a({\underline c})}(X*Y)\mod {\mathcal C}_\zeta^{< {\underline c}}$ 
where ${}^pH^i$ stands for perverse 
cohomology of degree $i$; here the associativity constraint comes
from the associativity of convolution and the unit object is 
$\bigoplus_{d\in {\mathcal D}({\underline c})} j_{d!*}^\zeta$ where
${\mathcal D}({\underline c})\subset {\underline c}$ is the subset 
of distinguished involutions, see~\cite{Lt}.

The category $\C_\zeta^{\underline c}$ contains a monoidal subcategory
$\C_\zeta^{\Gamma_{\underline c}\cap\Gamma_{\underline c}^{-1}}$ formed
by the direct sums of $j^\zeta_{w!*},\ w\in\Gamma_{\underline c}\cap
\Gamma_{\underline c}^{-1}$, see {\em loc. cit.}
Let $F_{\underline c}$ be the quotient of the centralizer of 
$u({\underline c})$ in $\LG_\zeta$ by its unipotent radical.
Note that the group $F_{\underline c}$ acts in a natural way on the set
$W_\zeta^0$.

The following Corollary of Proposition~\ref{key} (cf. Remark~\ref{some})
is proved similarly to Theorem 3 of~\cite{B}:

\begin{Cor} \label{corr}
There exists an equivalence of monoidal categories
$Coh_{F_{\underline c}}(W_\zeta^0\times W_\zeta^0)\simeq
\C_\zeta^{\Gamma_{\underline c}\cap\Gamma_{\underline c}^{-1}}$.
\end{Cor} \epf

We refer the reader to \cite{BO} 4.2 and 5.1 for the definition
of a $F_{\underline c}-$set $Y$ of centrally extended points and
of the monoidal category ${\mbox Fun}_{F_{\underline c}}(Y,Y)$.

\begin{Thm} 
\label{Jcat}
There exists a finite $F_{\underline c}-$set of centrally
extended points ${\bf X}$ and an equivalence of monoidal categories
$F: {\mathcal C}_\zeta^{\underline c}\to 
{\mbox Fun}_{F_{\underline c}}({\bf X}\times W_\zeta^0,{\bf X}\times W_\zeta^0)$.
Here $F_{\underline c}$ acts on ${\bf X}\times W_\zeta^0$ diagonally.
\end{Thm}

First let us consider the Serre subcategory ${}^0{\mathcal C}_\zeta^{\underline c}
\subset {\mathcal C}_\zeta^{\underline c}$ with simple objects of the form
$j_{w!*}^\zeta,\ w\in W(\LG^0_\zeta )\subset Stab_W(\zeta)$.
Let ${}^0F_{\underline c}\subset \LG_\zeta^0$ be the intersection of $F_{\underline c}$ 
and $\LG_\zeta^0$.
The arguments parallel to the proof of Theorem 4 in \cite{BO}
show that there exists a finite ${}^0F_{\underline c}-$set $\tilde Y$ of centrally extended points
and  a monoidal equivalence ${}^0{\mathcal C}_\zeta^{\underline c}\simeq Fun_{{}^0F_{\underline c}}
(\tilde Y, \tilde Y)$ (the proof involves the results parallel to 
Theorems 1, 2 and 3 from \cite{B}; the proofs are parallel to those
in \cite{B} using our Theorem \ref{main} instead of~\cite{KGB} used
in \cite{B}).

In particular the category ${}^0{\mathcal C}_\zeta^{\underline c}$ is rigid.
The category ${\mathcal C}_\zeta^{\underline c}$ is generated by 
${}^0{\mathcal C}_\zeta^{\underline c}$ and some invertible objects and hence is
rigid as well.

Note that the category ${\mathcal C}_\zeta^{\underline c}$ is not tensor but multi-tensor in
the sense that the unit object is decomposable (see e.g. \cite{ENO} 2.4 for a discussion of
the multi-fusion categories). It is easy to see that the category ${\mathcal C}_\zeta^{\underline c}$
is indecomposable, see {\em loc. cit.} The category 
$\C_\zeta^{\Gamma_{\underline c}\cap\Gamma_{\underline c}^{-1}}$ is a component category
of ${\mathcal C}_\zeta^{\underline c}$, which means that it is of the form 
$e\otimes {\mathcal C}_\zeta^{\underline c}\otimes e$ where $e$ is some direct summand 
of the unit object of the category ${\mathcal C}_\zeta^{\underline c}$. 
Thus the category ${\mathcal C}_\zeta^{\underline c}$ is dual to the category
$\C_\zeta^{\Gamma_{\underline c}\cap\Gamma_{\underline c}^{-1}}$ via the module
category ${\mathcal C}_\zeta^{\underline c}\otimes e$. In turn by Corollary \ref{corr} the
category $\C_\zeta^{\Gamma_{\underline c}\cap\Gamma_{\underline c}^{-1}}=
Coh_{F_{\underline c}}(W_\zeta^0\times W_\zeta^0)$ is dual to the category 
$Rep(F_{\underline c})$ via module category $Coh(W_\zeta^0)$. Hence the category
${\mathcal C}_\zeta^{\underline c}$ is dual to the category $Rep(F_{\underline c})$ via
some module category and therefore is of the form $Fun_{F_{\underline c}}(\tilde {\bf X},
\tilde {\bf X})$ for some finite $F_{\underline c}-$set $\tilde {\bf X}$ of centrally extended
points.

Now we show that actually $\tilde {\bf X}={\bf X}\times W_\zeta^0$ for some finite 
$F_{\underline c}-$set ${\bf X}$ of centrally extended points. Note first that $W_\zeta^0$
acts on the set $\tilde {\bf X}$. This follows from the fact that the category
${\mathcal C}_\zeta^{\underline c}$ contains invertible objects parametrized 
by the group $W_\zeta^0$ and on the other hand the category 
$Fun_{F_{\underline c}}(\tilde {\bf X},\tilde {\bf X})$ acts on $Coh(\tilde X)$ where 
$\tilde X$ is the finite set underlying $\tilde {\bf X}$. We claim that this action is
actually free. Indeed consider the fixed set $\tilde {\bf X}^w$ for some $1\ne w\in W_\zeta^0$.
This set is clearly $F_{\underline c}-$invariant and hence carries some $F_{\underline c}-$
equivariant sheaf $Y \in Coh_{F_{\underline c}}(\tilde {\bf X})$. Any irreducible constituent
of $\underline{\mbox{Hom}}(Y,Y)$  (see \cite{BO} 4.1) is clearly $w-$invariant 
(with respect to the right tensoring).
  But this is impossible since the category ${\mathcal C}_\zeta^{\underline c}$ is
  $W_\zeta^0-$graded and multiplication by $w$ just shifts the grading. Thus the set
$\tilde {\bf X}^w$ is empty and $W_\zeta^0-$action on $\tilde {\bf X}$ is free and hence
$\tilde {\bf X}={\bf X}\times W_\zeta^0$. Theorem is proved.
\epf

\subsection{Remark} Note that the set of orbits of $F_{\underline c}$ on the set 
${\bf X}\times W_\zeta^0$ identifies with the set of irreducible direct summands
of the unit object in ${\mbox Fun}_{F_{\underline c}}({\bf X}\times W_\zeta^0,{\bf X}\times W_\zeta^0)=
{\mathcal C}_\zeta^{\underline c}$ which in turn identifies with the set of distinguished
involutions in ${\underline c}$.


\subsection{Example} 
\label{43}
Let $G=SL_2$ and let $\zeta =\left( \begin{array}{cc}1&0\\ 0&-1
\end{array}\right)\in \LG=PGL_2$ be such that the
group $\LG_\zeta$ is disconnected. Then the algebra ${\mathcal H}_\zeta$ is isomorphic 
to the group algebra of the group
$\BZ\rtimes \BZ/2\BZ$ with the basis $C_w$ corresponding to the
group elements in the group algebra. In particular there is only one
cell ${\underline c}$ (left, right and two-sided) which coincides with 
the whole group.
The group $\LG_\zeta =F_{{\underline c}}$ is isomorphic 
to the semidirect product $\Gm\rtimes \BZ/2\BZ$. It is easy to see that in this
case the set ${\bf X}$ consists of two points (not centrally extended)
permuted transitively by the group $F_{\underline c}$.

\end{section}

\begin{section}{Cells in finite Weyl groups}

\subsection{Finite Weyl group} In this section we apply the previous results
to the truncated convolution category of monodromic sheaves on the finite
dimensional flag variety. The results are new even in the case when $\zeta$
is trivial.

The obvious embedding $G\subset G(F)$ induces embedding $G/B\subset \Fl$ 
and the torsor $\tFl$ restricts to the torsor $G/U\subset \tFl$. We consider
the full subcategory ${\mathcal C}_\zeta^f$ of the category 
${\mathcal C}_\zeta$ consisting of sheaves supported on $G/U$. 
The Grothendieck ring ${\mathcal K}_\zeta^f$ of the category 
${\mathcal C}_\zeta^f$ is identified with the subring of 
${\mathcal K}_\zeta\simeq{\mathcal H}_\zeta$ with basis $C_w$, 
$w\in W_f\cap Stab_W(\zeta)=:W_f(\zeta)$. Note that $W_f(\zeta)$ is
the semidirect product of the Coxeter group $W_f^\zeta$ (see~\ref{til})
with $W_\zeta^0$: we have $W_f(\zeta)=W_f^\zeta\rtimes W_\zeta^0$.
The notions of cells, $a$-function etc.
are defined for ${\mathcal H}_\zeta^f$ 
in the same way as for ${\mathcal H}_\zeta$.
In particular, for any two-sided cell ${\underline c}_f\subset W_f(\zeta)$
one defines the monoidal category ${\mathcal C}_\zeta^{{\underline c}_f}$
with truncated convolution as a tensor product in the same way as before.
The category ${\mathcal C}_\zeta^{{\underline c}_f}$ is multi-fusion category
in the sense of \cite{ENO}.

The cell ${\underline c}_f$ is contained in a unique two-sided cell
${\underline c}\subset Stab_W(\zeta)$. Let $A_{{\underline c}_f}$ be
the component group of the group $F_{{\underline c}}$. We have the
following

\begin{Prop} \label{grth}
The multi-fusion category ${\mathcal C}_\zeta^{{\underline c}_f}$
is group-theoretical. Moreover there exists a subquotient $S$ of the group $A_{{\underline c}_f}$
and 3-cocycle $\omega \in Z^3(S,\Gm)$ such that the category ${\mathcal C}_\zeta^{{\underline c}_f}$
is dual to the category $Vec_{S,\omega}$ with respect to a module category (here 
$Vec_{S,\omega}$ is the category of coherent sheaves on $S$ with convolution tensor
product and associativity defined by $\omega$, see \cite{ENO} 8.8).
\end{Prop}

{\bf Proof} is similar to the proof of Proposition 8.44(i) in \cite{ENO} since the category
${\mathcal C}_\zeta^{{\underline c}_f}$ is a tensor subcategory of the category
${\mathcal C}_\zeta^{{\underline c}}$. 
\epf

\subsection{Remark}
Recall that we have a map $A_{{\underline c}_f}\to W^0_\zeta$.
It is easy to see
that the subquotient $S$ is such that the map $S\to W^0_\zeta$
is well defined. Moreover, the
image of $S\to W^0_\zeta$ coincides with the
image of $A_{{\underline c}_f}\to W^0_\zeta$.

\subsection{}
Now recall that for a two-sided cell 
${\underline c}_f\subset W_f(\zeta)$ Lusztig defined a finite group
$\G_{{\underline c}_f}$ together with homomorphism $\G_{{\underline c}_f}\to W^0_\zeta$, 
see \cite{l107} for the case $W^0_\zeta=1$ and \cite{l1} IV for the general case. We expect 
that in Proposition \ref{grth} we can always choose $S=\G_{{\underline c}_f}$; moreover $\omega$ should be trivial 
except, possibly, the case when ${\underline c}_f$ is associated with exceptional family. Recall
that there are just 2 exceptional families
for $W_f^\zeta$ of type $E_8$, and one more for $W_f^\zeta$ of type $E_7$.

\begin{Conj}
\label{finite}
Let ${\underline c}_f\subset W_f(\zeta)$ be 
a two-sided cell such that the corresponding family is not exceptional. Then the category
 ${\mathcal C}_\zeta^{{\underline c}_f}$ is tensor equivalent to
$Fun_{\G_{{\underline c}_f}} 
({\bf X}_f\times W^0_\zeta, {\bf X}_f\times W^0_\zeta)$
for some finite $\G_{{\underline c}_f}-$set ${\bf X}_f$ of centrally extended points. 
\end{Conj}

In support of Conjecture~\ref{finite} we can state

\begin{Thm}
\label{finit}
Conjecture~\ref{finite} is true in the following cases:

(a) $W^0_\zeta =1$;

(b) the map $\G_{{\underline c}_f}\to W^0_\zeta$ is injective.
\end{Thm}

\proof (a):
It is enough to prove that the category
${\mathcal C}_\zeta^{{\underline c}_f}$ is dual to the category $Vec_{\G_{{\underline c}_f}}=Coh(\G_{{\underline c}_f})$ (with
tensor product given by convolution and obvious associativity constraint)
with respect to some module category.

For a subset $A\subset {\underline c}_f$ let ${\mathcal C}_\zeta^A$ 
denote the subcategory
(usually not tensor) of ${\mathcal C}_\zeta^{{\underline c}_f}$ 
formed by the direct sums
of $j^\zeta_{w!*},\ w\in A$. 
Let $\Gamma \subset {\underline c}_f$ be a left cell. Then the
category ${\mathcal C}_\zeta^{\Gamma \cap \Gamma^{-1}}$ is 
a tensor subcategory of
${\mathcal C}_\zeta^{{\underline c}_f}$; moreover it is a component category of
${\mathcal C}_\zeta^{{\underline c}_f}$, see \cite{ENO} 2.4. 
To prove the theorem it is enough
to show that some
component category of ${\mathcal C}_\zeta^{{\underline c}_f}$ is dual
to the category $Vec_{\G_{{\underline c}_f}}$, see {\em loc. cit.} 5.5.

Recall that for any left cell $\Gamma \subset W_f(\zeta)$
G.~Lusztig defined the associated representation
$[\Gamma]$ of $W_f(\zeta)$, see \cite{l107}. 
All possible representations of the form
$[\Gamma]$ were computed by Lusztig in \cite{lfrench}. 
We are going to use these results
together with the following fact (see \cite{l107} 12.15):

(*) for two left cells $\Gamma_1$ and $\Gamma_2$ 
the cardinality of intersection
$\Gamma_1\cap \Gamma_2^{-1}$ equals $\dim Hom ([\Gamma_1],[\Gamma_2])$.

Recall also that the Frobenius-Perron dimension (see \cite{ENO}) 
of a fusion category 
${\mathcal C}$ can be read of the character table of the Grothendieck ring of ${\mathcal C}$.
  For the categories ${\mathcal C}_\zeta^{\Gamma \cap \Gamma^{-1}}$ these character
  tables are known, see~\cite{llead},~3.14. 
We see that in all cases the Frobenius-Perron
  dimension of ${\mathcal C}_\zeta^{\Gamma \cap \Gamma^{-1}}$ 
is equal to $|\G_{{\underline c}_f}|$. Thus
  for the group $S$ from Proposition~\ref{grth}
we have $|S|=|\G_{{\underline c}_f}|$.

We now proceed case by case. Assume first that the group $\G_{{\underline c}_f}$ is elementary abelian 2-group
(possibly trivial), but ${\underline c}_f$ is not exceptional. It is known (see \cite{llead} 3.11) that 
for any left cell $\Gamma$ the category ${\mathcal C}_\zeta^{\Gamma \cap \Gamma^{-1}}$
is tensor equivalent to $Vec_{S,\omega}$ where the group $S$ is isomorphic (non canonically)
to $\G_{{\underline c}_f}$. We just need to show that $\omega$ is trivial. For this let us note that the results
of \cite{lclass} and \cite{lfrench} imply that we always have another left cell $\Gamma_1$
such that the cardinality $\Gamma \cap \Gamma_1^{-1}$ is 1 (this is exactly what fails for
exceptional cells). But obviously ${\mathcal C}_\zeta^{\Gamma \cap \Gamma_1^{-1}}$
is a module category over ${\mathcal C}_\zeta^{\Gamma \cap \Gamma^{-1}}$ and hence
${\mathcal C}_\zeta^{\Gamma \cap \Gamma^{-1}}=Vec_{S,\omega}$ admits a tensor functor
to the category of vector spaces. This forces $\omega$ to be trivial 
and the theorem is proved in this case.

Now assume that $\G_{{\underline c}_f}$ is symmetric group in 3 letters $S_3$. The results
of \cite{lclass} and \cite{lfrench} imply that we always have two left cells $\Gamma_1$
and $\Gamma_2$ such that $|\Gamma_1\cap \Gamma_1^{-1}|=|\Gamma_2\cap \Gamma_2^{-1}|=3$
and $|\Gamma_1\cap \Gamma_2^{-1}|=2$. The existence of $\Gamma_1$ implies that the
group $S$ from  Proposition~\ref{grth} 
is isomorphic to $S_3$ (since any category
dual to $Vec_{\BZ/6\BZ ,\omega}$ has 6 simple objects) and we just need to show
that $\omega$ is trivial. The category ${\mathcal C}_\zeta^{\Gamma_1 \cap \Gamma_1^{-1}}$
is of the form $\C (S_3, H, \omega ,\psi)$
(see~\cite{ENO},~8.8) for some subgroup $H\subset S_3$,
3-cocycle $\omega \in Z^3(S_3,\Gm)$ and 2-cochain $\psi \in C^2(H,\Gm)$ such that
$\omega|_H=d(\psi)$. It is easy to see that
either $H=S_3$ or $H=\BZ/2\BZ$ (otherwise we will have 6 simple objects in the category
$\C (S_3, H, \omega ,\psi)$). In the first case $\omega$ is trivial and we are done;
in the second case consider the category ${\mathcal C}_\zeta^{\Gamma_2 \cap \Gamma_2^{-1}}$.
It is also of the form $\C (S_3, H', \omega ,\psi')$ and $H'\ne \BZ/2\BZ$ since otherwise
$\Gamma_1\cap \Gamma_2$ would have 3 simple objects. Hence $H'=S_3$ and the
theorem is proved in this case.

Next assume that $\G_{{\underline c}_f} =S_4$. 
In this case $\zeta =1$, $\LG$ is of type $F_4$, and 
$A_{\underline c}=S_4$. Since $S$ is a
subquotient of $A_{\underline c}$ such that
$|S|=24$, we have $S=S_4$. Let $\Gamma_1, \Gamma_2$ be left cells such that
$[\Gamma_1]=\chi_{12,1}+\chi_{16,1}+\chi_{9,2}+\chi_{6,1}+\chi_{4,3}$ and
$[\Gamma_2]=\chi_{12,1}+\chi_{16,1}+\chi_{9,3}+\chi_{6,1}+\chi_{4,4}$
(notations of~\cite{lclass}). 
Then $|\Gamma_1\cap \Gamma_1^{-1}|=|\Gamma_2\cap \Gamma_2^{-1}|=5$
and $|\Gamma_1\cap \Gamma_2^{-1}|=3$. The character table of the Grothendieck ring of the
category  ${\mathcal C}_\zeta^{\Gamma_1 \cap \Gamma_1^{-1}}$ is contained in \cite{llead}
and coincides with the character table of the Grothendieck ring of the
category  ${\mathcal C}_\zeta^{\Gamma_2 \cap \Gamma_2^{-1}}$.
This character table uniquely determines the structure of this Grothendieck ring as a
based ring since it is commutative; we see that this Grothendieck ring coincides with
the Grothendieck ring of the category of representations of $S_4$. In particular it
contains just two nontrivial based subrings: representations of $S_2$ and representations
of $S_3$ of Frobenius-Perron dimensions 2 and 6 respectively. Now the category
${\mathcal C}_\zeta^{\Gamma_1 \cap \Gamma_1^{-1}}$ is of the form
$\C (S_4, H, \omega ,\psi)$; hence it contains subcategory $Rep(H)$ of Frobenius-Perron
dimension $|H|$. Thus $|H|=1$, or $|H|=2$, or $|H|=6$, or $|H|=24$. In the first
two cases the category $\C (S_4, H, \omega ,\psi)$ contains at least $24/4=6>5$ objects
and this is impossible. In the last case $H=S_4$ and therefore $\omega$ is trivial. 
Assume that $|H|=6$ and therefore $H=S_3\subset S_4$.
Consider the category ${\mathcal C}_\zeta^{\Gamma_2 \cap \Gamma_2^{-1}}$.
It is also of the form $\C (S_4, H', \omega ,\psi')$ and $H'\ne S_3$ since otherwise
$\Gamma_1\cap \Gamma_2^{-1}$ would have 5 simple objects. Hence $H'=S_4$ and 
the theorem is proved in this case.

Finally assume that $\G_{{\underline c}_f} =S_5$. In this case $\zeta =1$, $\LG$ is of type $E_8$, and 
$A_{{\underline c}_f}=S_5$. 
Since $S$ is subquotient of $A_{{\underline c}_f}$ such that
$|S|=120$, we have $S=S_5$. Let $\Gamma$ be a left cell such that
$[\Gamma]=4480_y+3150_y+4200_y+420_y+7168_w+1344_w+2016_w$ (notations
of~\cite{lclass}). The character table of the Grothendieck ring of the
category  ${\mathcal C}_\zeta^{\Gamma \cap \Gamma^{-1}}$ is contained 
in~\cite{llead}. 
This character table uniquely determines the structure of this Grothendieck ring as a
based ring since it is commutative; we see that this Grothendieck ring coincides with
the Grothendieck ring of the category of representations of $S_5$ and in particular it
contains just one nontrivial based subring of Frobenius-Perron dimension 2. On the
other hand the category ${\mathcal C}_\zeta^{\Gamma \cap \Gamma^{-1}}$ is of
the form $\C (S_5, H, \omega ,\psi)$ and hence contains subcategory $Rep(H)$ of
Frobenius-Perron dimension $|H|$. Thus either $|H|=1$, or $|H|=2$, or $|H|=120$.
In the first two cases the category $\C (S_5, H, \omega ,\psi)$ has at least $120/4=30>7$
simple objects and this is impossible. Thus $H=S_5$ and hence $\omega$ is trivial.
The theorem is proved in this case and (a) is proved.

(b) follows easily from Theorem \ref{Jcat} since in this case the category
${\mathcal C}_\zeta^{\Gamma \cap \Gamma^{-1}}$ is equivalent to
${\mathcal W}_0^\zeta$.

\epf

\subsection{Remark} Part (a) of the Theorem covers all the 
cases when group $G$ is adjoint, or, more
generally, has connected center. Part (b) covers all the cases when the group $G$ is isogenous
to a product of groups of type $A$.

\subsection{Example} 
Let $G$ and $\zeta$ be the same as in Example~\ref{43}.
Then the category ${\mathcal C}_\zeta^{{\underline c}_f}$
contains two simple objects ${\bf 1}$ and $\delta$ where ${\bf 1}$ is the
unit object and $\delta \bullet \delta ={\bf 1}$. It is well known that
there are two nonequivalent monoidal categories of this kind: one is
$Rep(\BZ/2\BZ)$ and the second differs from the first by the twist of the
associativity constraint by a 3-cocycle. Theorem \ref{finit} asserts in
this case that the category  ${\mathcal C}_\zeta^{{\underline c}_f}$
is equivalent to $Rep(\BZ/2\BZ)$.

\subsection{Lusztig's Conjecture}
 Let us assume that $W^0_\zeta$ is trivial. G.~Lusztig conjectured (see \cite{llead} 3.15) that
for any left cell $\Gamma \subset {\underline c}_f$ there is an isomorphism of based rings
$K({\mathcal C}_\zeta^{\Gamma \cap \Gamma^{-1}})\simeq K_{\G_{{\underline c}_f}}
(\G_{{\underline c}_f}/\H \times \G_{{\underline c}_f} /\H )$ for a suitable subgroup
$\H \subset \G_{{\underline c}_f}$ (list of possible subgroups $\H$ is contained in {\em loc. cit.}).
Moreover, this Conjecture is verified by Lusztig in {\em loc. cit.}
for the case when $\G_{{\underline c}_f}$ is elementary abelian and in general (unpublished). 
We would like to point out how our results imply this Conjecture.
Actually, we will prove more: in almost all cases there exists an equivalence of tensor categories 
${\mathcal C}_\zeta^{\Gamma \cap \Gamma^{-1}}\simeq Coh_{\G_{{\underline c}_f}}
(\G_{{\underline c}_f}/\H \times \G_{{\underline c}_f} /\H )$ (this was conjectured by
Lusztig in \cite{Lt}). We are going to prove the following even more general statement
(also conjectured by Lusztig):

\begin{Thm} \label{lcon}
Assume that $W^0_\zeta$ is trivial.

(a) (cf. \cite{llead}) For any two-sided cell ${\underline c}_f\subset W_f(\zeta)$ there exists a finite
$\G_{{\underline c}_f}-$set $X$ (with no central extensions involved!) and an isomorphism
of based rings $K({\mathcal C}_\zeta^{{\underline c}_f})\simeq 
K_{\G_{{\underline c}_f}}(X\times X)$. 

(b) (cf. \cite{Lt}) Assume in addition that the family attached to ${\underline c}_f$ is not exceptional.
Then there exists an equivalence of tensor categories
${\mathcal C}_\zeta^{{\underline c}_f}\simeq Coh_{\G_{{\underline c}_f}}(X\times X)$. 
\end{Thm}

\proof 
First of all note that part (b) implies part (a) since for exceptional cells ${\underline c}_f$
we have $\G_{{\underline c}_f}=\BZ/2\BZ$ and the result is easy to deduce. Thus we
will prove only (b). 

It follows from Theorem \ref{finit} (a) that there exists a module category
$\M$ over $Rep(\G_{{\underline c}_f})$ such that the category ${\mathcal C}_\zeta^{{\underline c}_f}$
is tensor equivalent to $Fun_{Rep(\G_{{\underline c}_f})}(\M,\M)$. Moreover, the left
cells in ${\underline c}_f$ are in bijection with indecomposable direct summands of $\M$
and for two indecomposable direct summands $\M_1$, $\M_2$ corresponding
to cells $\Gamma_1$, $\Gamma_2$ we have $|\Gamma_1\cap \Gamma_2^{-1}|=$ number of
simple objects in the category $Fun_{Rep(\G_{{\underline c}_f})}(\M_1,\M_2)$. 
Recall that any
indecomposable module category over $Rep(\G_{{\underline c}_f})$ is of the form
$Rep^1(\tilde H)$ where $H\subset \G_{{\underline c}_f}$ is a subgroup and $\tilde H$ is
a central extension of $H$ by $\Gm$, see \cite{BO}. We just need to show that we can choose
$\M$ in such a way that for any indecomposable direct summand $\M_1$ of $\M$ the
corresponding central extension $\tilde H$ splits. Now we will proceed case by case
(in the cases $\G_{{\underline c}_f}=S_5$ or $\G_{{\underline c}_f}=S_4$ we are going to use
the tables of module categories over $Rep(\G_{{\underline c}_f})$, see Appendix).

Let us consider first the case $\G_{{\underline c}_f}=S_5$.
Let $\Gamma_0\subset {\underline c}_f$ be the left cell
from the proof of Theorem \ref{finit} (that is 
$[\Gamma_0]=4480_y+3150_y+4200_y+420_y+7168_w+1344_w+2016_w$).
We have seen in the proof of Theorem \ref{finit} that the category
${\mathcal C}_\zeta^{\Gamma_0 \cap \Gamma_0^{-1}}$ is equivalent to 
$Rep(S_5)$. It is known (see \cite{ENO}) that we can choose $\M =
{\mathcal C}_\zeta^{\Gamma_0}$ (this is considered as a module category
over ${\mathcal C}_\zeta^{\Gamma_0 \cap \Gamma_0^{-1}}$). Let $\Gamma_1\subset
{\underline c}_f$ be any left cell, let $\M_1$ be the corresponding indecomposable
summand of $\M$, let $H\subset S_5$ be
the corresponding subgroup and assume that
the central extension $\tilde H$ is described by 2-cocycle $\psi \in H^2(H,\Gm)$.
Note that the number of simple objects in the category $Fun_{Rep(\G_{{\underline c}_f})}(\M_1,\M_1)$=
number of simple objects in $\C (S_5, H, 1 ,\psi)$. It is easy to see that for any $H\subset S_5$ 
the cardinality of $H^2(H,\Gm)$ is $\le 2$; this implies that the number of simple objects in 
$\C (S_5, H, 1 ,\psi)$ is the same as in $\C (S_5, H, 1 ,1)$. Hence
$|\Gamma_1 \cap \Gamma_1^{-1}|=$number of simple objects in $\C (S_5, H, 1 ,1)$.
Also $|\Gamma_1 \cap \Gamma_0^{-1}|=$ number of
simple objects in the category $Fun_{Rep(S_5)}(Rep(S_5),\M_1)=$
number of irreducible projective representations of $H$ with
respect to cocycle $\psi$. An explicit classification of subgroups $H\subset S_5$ (there are 19
such subgroups up to conjugacy) together with classification of possible modules of the form
$[\Gamma]$ from \cite{lfrench} shows that $H$ should be from the list \cite{llead} 2.13 and 
$\psi$ should be trivial. This case is finished.

Now consider the case $\G_{{\underline c}_f}=S_4$. Let $\Gamma_0$ be a left cell such that
$[\Gamma_0]=\chi_{12,1}+\chi_{16,1}+\chi_{9,2}+\chi_{6,1}+\chi_{4,3}.$
In the proof of
Theorem \ref{finit} we showed that the category ${\mathcal C}_\zeta^{\Gamma_0 \cap \Gamma_0^{-1}}$
is equivalent to $\C=Rep(S_4)$ (we showed that it is equivalent either to $Rep(S_4)$ or to
$\C (S_4, S_3, 1 ,1)$ but these two categories are in fact equivalent). As before we take
$\M ={\mathcal C}_\zeta^{\Gamma_0}$. Now for any left cell $\Gamma_1$ such that
$[\Gamma_1]=\chi_{12,1}+\chi_{16,1}+\chi_{9,2}+\chi_{6,1}+\chi_{4,3}$ the corresponding
indecomposable module category $\M_1$ is of the form $Rep(S_4)$ (this is the only indecomposable 
module category $\M_1$ over $\C$ with 5 irreducible objects and such that $Fun_\C(\M_1,\M_1)$
has 5 objects). Similarly let $\Gamma_2$ be a cell such that 
$[\Gamma_2]=\chi_{12,1}+\chi_{16,1}+2\chi_{9,2}+\chi_{6,2}+\chi_{4,3}+\chi_{1,2}$;
then the corresponding module category $\M_2$ over $\C$ is of the form 
$Rep(D_8)$ where $D_8$ is a subgroup of order 8 (this is the only indecomposable 
module category $\M_2$ over $\C$ with 5 irreducible objects and such that $Fun_\C(\M_2,\M_2)$
has 9 objects). Now let $\Gamma_3$ be a {\em unique} cell (see \cite{l107} 12.12) such that 
$[\Gamma_3]=\chi_{12,1}+\chi_{16,1}+2\chi_{9,3}+\chi_{6,2}+\chi_{4,4}+\chi_{1,3}$;
then the corresponding module category $\M_3$ over $\C$ has 2 irreducible objects and 
the category $Fun_\C(\M_3,\M_3)$ has 9 irreducible objects. Thus either $\M_3$ is of the
form  $Rep(S_2)$ where $S_2$ is the symmetric group in two letters or is of the form
$Rep^1(\tilde D_8)$ where $\tilde D_8$ is the nontrivial central extension of $D_8$. One
observes that these two module categories are interchanged by a tensor autoequivalence
of $Rep(S_4)$ (and any autoequivalence fixes the module categories $\M_1$ and $\M_2$)
\footnote{Let $Kl$ be the Klein's subgroup of $S_4$ and let $\widetilde{Kl}$ be its nontrivial
central extension. The module category $\N =Rep^1(\widetilde{Kl})$ has one simple object
and the dual category $\C^*_{\N}$ is pointed, i.e. all its simple objects are invertible. Moreover,
the group of isomorphism classes of invertible objects in $\C^*_{\N}$ is isomorphic to $S_4$.
A {\em choice} of such an isomorphism produces the desired autoequivalence.}.
Thus we can assume that $\M_3$ is of the form $Rep(S_2)$.
 Next let $\Gamma_4$ be a cell such that 
$[\Gamma_4]=\chi_{12,1}+\chi_{16,1}+\chi_{9,3}+\chi_{6,1}+\chi_{4,4}$;
then the corresponding module category $\M_4$ over $\C$ is of the form 
$Rep(S_3)$ where $S_3$ is the symmetric group in three letters (this is the only indecomposable 
module category $\M_4$ over $\C$ with 3 irreducible objects, such that $Fun_\C(\M_4,\M_4)$
has 5 objects and $Fun(\M_3,\M_4)$ has 5 objects). Finally, let $\Gamma_5$ be a cell such that 
$[\Gamma_5]=\chi_{12,1}+2\chi_{16,1}+\chi_{9,2}+\chi_{9,3}+\chi_{6,2}+\chi_{4,1}$;
then the corresponding module category $\M_5$ over $\C$ is of the form 
$Rep(S_2\times S_2)$ (this is the only indecomposable 
module category $\M_5$ over $\C$ with 4 irreducible objects, such that $Fun_\C(\M_5,\M_5)$
has 9 objects and $Fun(\M_4,\M_5)$ has 4 objects). It is proved in \cite{l107} 12.12 that any left
cell from ${\underline c}_f$ was listed above and we are done in this case.

The case $\G_{{\underline c}_f}=S_3$ is trivial since for any subgroup $H\subset S_3$ we
have $H^2(H,\Gm)=0$.

Now let us consider the case when $\G_{{\underline c}_f}$ is elementary abelian group.
To any left cell $\Gamma \subset {\underline c}_f$ Lusztig associated a subgroup
$\H_\Gamma \subset \G_{{\underline c}_f}$ such that for any two cells $\Gamma_1$,
$\Gamma_2$ we have $|\Gamma_1\cap \Gamma_2^{-1}|=$ number of simple objects
in the category $Coh_{\G_{{\underline c}_f}}(\G_{{\underline c}_f}/\H_{\Gamma_1}\times
\G_{{\underline c}_f}/\H_{\Gamma_2})$, see \cite{llead} 3.16 (a). It is known that
there exist cells $\Gamma_0$ and $\Gamma_1$ such that
$\H_{\Gamma_0}=\G_{{\underline c}_f}$ and $\H_{\Gamma_1}=\{ 0\}$ (recall that we assume
that the family associated with ${\underline c}_f$ is not exceptional). Thus the category
${\mathcal C}_\zeta^{\Gamma_0 \cap \Gamma_0^{-1}}$ can be identified with the
category $Rep(\G_{{\underline c}_f})$. The module category 
${\mathcal C}_\zeta^{\Gamma_0 \cap \Gamma_1^{-1}}$ over 
${\mathcal C}_\zeta^{\Gamma_0 \cap \Gamma_0^{-1}}$ has
one irreducible object; we can assume
that it is of the form $Rep^1(\tilde H)$ where $\tilde H$ is the trivial central extension of the
trivial group $H=\H_{\Gamma_1}=\{ 0\}$ (indeed, it is easy to see that the action the group 
of tensor autoequivalences of the category $Rep(\G_{{\underline c}_f})$ on the set
of module categories with one irreducible object is transitive). Now let $\Gamma_2 \subset 
{\underline c}_f$ be a left cell, let $\M_2$ be the corresponding indecomposable
summand of $\M$, let $H(\Gamma_2)\subset \G_{{\underline c}_f}$ be
the corresponding subgroup
and let $\tilde H(\Gamma_2)$ be the corresponding central extension. We have
$|\H_{\Gamma_2}|=|\Gamma_0 \cap \Gamma_2^{-1}|=$ number of irreducible objects
in $Rep^1(\tilde H(\Gamma_2))$, and $\frac{|\G_{{\underline c}_f}|}{|\H_{\Gamma_2}|}=
|\Gamma_1 \cap \Gamma_2^{-1}|=\frac{|\G_{{\underline c}_f}|}{|H(\Gamma_2)|}$ whence the
number of irreducible objects in $Rep^1(\tilde H(\Gamma_2))$ equals to $|\H_{\Gamma_2}|$
equals to $|H(\Gamma_2)|$. This obviously implies that the central extension
$\tilde H(\Gamma_2)$ splits and we are done (note also that now \cite{llead} Proposition 3.8
implies that $H(\Gamma)=\H_\Gamma$
for any left cell $\Gamma \subset {\underline c}_f$). \epf

\subsection{Remark} Let ${\underline c}_f$ be a two-sided cell such that associated family is
exceptional. The results of \cite{l107} (see e.g. Theorem 11.2) 
suggest that for any left cell $\Gamma \subset
{\underline c}_f$ the category  ${\mathcal C}_\zeta^{\Gamma \cap \Gamma^{-1}}$ is
equivalent to the category $Coh(\BZ/2\BZ)$ with convolution tensor product and
associativity constraint given by the nontrivial 3-cocycle.

\end{section}

\begin{section}{Applications to character sheaves}

\subsection{} Let $\tD^\zeta_f$ stand for the full monoidal subcategory
of $\tD^\zeta_I$ formed by the complexes supported at $G/U\subset\tFl$.
Recall the horocycle space $\CY$: the quotient of 
$G/U\times G/U$ modulo the diagonal right action of $T$ (see~\cite{MV}).
It is a $G$-equivariant $T$-torsor over $G/B\times G/B$, and 
$\tD^\zeta_f$ is equivalent to the monoidal (with respect to convolution)
$G$-equivariant constructible derived category of $\zeta$-monodromic
complexes on $\CY$. So from now on we will consider $\tD^\zeta_f$ in the
latter incarnation. Let $D_G(G)$ stand for the $G$-equivariant 
(with respect to the adjoint action) constructible derived category on $G$.
Recall the averaging functor $\Gamma^G:\ \tD^\zeta_f\to D_G(G)$ of~\cite{MV}.
By definition, the character sheaves on $G$ with central character $\zeta$
are the irreducible constituents of the perverse cohomology sheaves
of the various objects of the form $\Gamma^G(\CF)$ where $\CF\in\tD_f^\zeta$,
see~\cite{MV}. These character sheaves are classified in~\cite{l1} in
terms of cells ${\underline{c}}_f\subset W_f(\zeta)$, and nonabelian
Fourier transform (see~\cite{l100}) attached to 
the finite group $\G_{{\underline c}_f}$. Equivalently, due to the
results of~\cite{llead}, this classification can be formulated in terms
of representations
of the Drinfeld double of the finite group $\G_{{\underline c}_f}$.
Our goal is to place the Lusztig's classification into the context of
finite monoidal categories, thus explaining the appearance of the Drinfeld
double.

Recall (see \cite{MV})  that the functor $\Gamma^G$ has the right 
adjoint functor $\Gamma^U$.
Recall that both categories 
$\tD^\zeta_f$ and $D_G(G)$ have natural monoidal structures 
(with respect to the
convolution $*$) and the functor $\Gamma^U$ has an obvious tensor structure. 
What is more important for us, 
the functor $\Gamma^U$ has an obvious structure of {\em central} functor, that 
is it factors through the
Drinfeld center $Z(\tD^\zeta_f)$. The following statement is elementary:

\begin{Prop} 
\label{nsf}
Let $\A$ be a category (with no monoidal structure) 
and let $\C$
be a pivotal rigid monoidal category. Let $(\fG,\fF)$ be a pair of adjoint functors 
$\fF:\A \to \C$ and 
$\fG:\C \to \A$. Then the structures of central functor on $\fF$ (that is 
factorizations through
$Z(\C)$) are in bijection with the following structures on the functor
$\fG$: the functorial
isomorphisms $u: \fG(X\otimes Y)\to \fG(Y\otimes X)$
such that $u$ is the identity
when $Y$ is the unit object of $\C$ and the composition 
$$\fG((X\otimes Y)\otimes Z)\to \fG(Z\otimes (X\otimes Y))=
\fG((Z\otimes X)\otimes Y)\to \fG(Y\otimes (Z\otimes X))$$ 
coincides with the morphism 
$$\fG((X\otimes Y)\otimes Z)=\fG(X\otimes (Y\otimes Z))\to \fG((Y\otimes Z)\otimes X)
=\fG(Y\otimes (Z\otimes X)).$$
\end{Prop}

We will call a functor $\fG$ from a monoidal category to an arbitrary category 
endowed 
with a structure described in Proposition~\ref{nsf} a {\em
commutator functor}. Thus the right adjoint of a commutator functor is a 
central functor and
the left adjoint of a central functor is a commutator functor. One shows
from the above
that the functor $\Gamma^G$ has a natural structure of the commutator 
functor.

Let $\C^\zeta_G(G)\subset D_G(G)$
be the semisimple abelian category formed by the direct sums of
character sheaves with the central character $\zeta$.
Obviously, for a finite cell ${\underline{c}}_f\subset W_f(\zeta)$, the
restriction of the functor $\Gamma^G$ to the category 
$\C_\zeta^{{\underline{c}}_f}\subset\C_\zeta^f\subset\tD_f^\zeta$
still has a structure of commutator functor. Moreover, it is easy 
to see from the definitions that the
functor ${}^pH^{a({\underline{c}}_f)}\Gamma^G$ from
$\C_\zeta^{{\underline{c}}_f}$ lands into $\C^\zeta_G(G)$, and has a 
structure
of the commutator functor with respect to the monoidal 
structure $\bullet$ on the
category $\C_\zeta^{{\underline{c}}_f}$. 
Recall that Lusztig \cite{l1} defined a 
decomposition
$\C^\zeta_G(G)=\bigoplus_{{\underline{c}}_f} 
\C^\zeta_G(G)^{{\underline{c}}_f}$ where the summands are also labeled by the 
two sided
cells in $W_f(\zeta)$.

Let $\fG: 
\C_\zeta^{{\underline{c}}_f}\to
\C^\zeta_G(G)^{{\underline{c}}_f}$ be the summand of the functor
${}^pH^{a({\underline{c}}_f)}\Gamma^G$
corresponding to the direct summand
$\C^\zeta_G(G)^{{\underline{c}}_f}\subset \C^\zeta_G(G)$. 
It follows from the
previous discussion that $\fG$ has a natural structure of the commutator 
functor with respect
to the monoidal structure $\bullet$. Thus its right adjoint functor 
$\tilde \fF:\ \C^\zeta_G(G)^{{\underline{c}}_f}\to
\C_\zeta^{{\underline{c}}_f}$ 
has a structure of the central functor (indeed it 
is known that the
monoidal structure $\bullet$ is rigid); or in other words the functor 
$\tilde \fF$ factors through
the functor $\fF: \C^\zeta_G(G)^{{\underline{c}}_f}\to 
Z(\C_\zeta^{{\underline{c}}_f})$ 
(here $Z(\C_\zeta^{{\underline{c}}_f})$ is the
Drinfeld center of the category $\C_\zeta^{{\underline{c}}_f}$).

\begin{Conj}
\label{vain} 
The functor $\fF$ is an equivalence of categories.
\end{Conj}

Now according to Conjecture~\ref{finite}, the monoidal category
$\C_\zeta^{{\underline{c}}_f}$ is monoidally equivalent to\\ 
$Fun_{\G_{{\underline c}_f}}({\bf X}_f\times W_\zeta^0,
{\bf X}_f\times W_\zeta^0)$.
Also, according to~\cite{O}, 
we have an equivalence of monoidal categories 
$Z(Fun_{\G_{{\underline c}_f}}({\bf X}_f\times W_\zeta^0,
{\bf X}_f\times W_\zeta^0))\simeq
Z(Rep(\G_{{\underline c}_f}))$.
Thus Conjecture~\ref{vain} implies that the simple objects in 
$\C^\zeta_G(G)^{{\underline{c}}_f}$ are labeled by
simple objects in $Z(Rep(\G_{{\underline c}_f}))$. 
This should be exactly Lusztig's parametrization 
of the character sheaves from~\cite{l1}.

\end{section}

\section{Appendix: Module categories over $Rep(S_4)$ and $Rep(S_5)$} 
For the reader's convenience we give here the tables of indecomposable module categories over
tensor categories $Rep(S_4)$ and $Rep(S_5)$. Any such category is of the form
$Rep^1(\tilde H)$ where $\tilde H$ is a central extension of subgroup $H$; in the case
when the central extension is trivial we have $Rep^1(\tilde H)=Rep(H)$; thus tilde in the
tables below always refer to the nontrivial extension (we have at most one nontrivial extension
in each case). The second column of the tables gives the number of irreducible objects
in $\M$; the third column gives the number of irreducible objects in $Fun(\M,\M)$. 

The notation $\langle g\rangle$ refers to the cyclic group generated by $g$; 
$S_i$ (or $S_i\times S_j$)
denotes to the standard parabolic subgroups; $Kl$ is the Klein's four-group; $H_{10}$ and $H_{20}$
are unique up to conjugacy subgroups of order 10 and 20; finally $H_6$ is a nonabelian subgroup
of order 6 in $S_3\times S_2$ which is not $S_3\times \{ e\}$.

\newpage

\centerline{Module categories over $Rep(S_4)$.}

$$\begin{array}{|c|c|c|}
\hline \M & \# \M &\# Fun(\M,\M) \\
\hline Rep(\{ e\}) & 1 & 24\\
\hline Rep(S_2)& 2 & 9\\
\hline Rep(\langle(12)(34)\rangle)& 2 & 12\\
\hline Rep(\langle(123)\rangle)& 3 & 8\\
\hline Rep(S_2\times S_2)& 4 & 9\\
\hline Rep^1(\widetilde{S_2\times S_2})& 1 & 9\\
\hline Rep(Kl)& 4 & 24\\
\hline Rep^1(\widetilde{Kl})& 1 & 24\\
\hline Rep(\langle(1234)\rangle)& 4 & 9\\
\hline Rep(S_3)& 3 & 5\\
\hline Rep(D_8)& 5 & 9\\
\hline Rep^1(\widetilde{D_8})& 2 & 9\\
\hline Rep(A_4)& 4 & 8\\
\hline Rep^1(\widetilde{A_4})& 3 & 8\\
\hline Rep(S_4)& 5 & 5\\
\hline Rep^1(\widetilde{S_4})& 3 & 5\\
\hline \end{array}$$

\newpage

\centerline{Module categories over $Rep(S_5)$.}

$$
\begin{array}{|c|c|c|}
\hline \M & \# \M &\# Fun(\M,\M) \\
\hline Rep(\{ e\}) & 1 & 120\\
\hline Rep(S_2)& 2 & 39\\
\hline Rep(\langle(12)(34)\rangle)& 2 & 36\\
\hline Rep(\langle(123)\rangle)& 3 & 24\\
\hline Rep(S_2\times S_2)& 4 & 21\\
\hline Rep^1(\widetilde{S_2\times S_2})& 1 & 21\\
\hline Rep(Kl)& 4 & 30\\
\hline Rep^1(\widetilde{Kl})& 1 & 30\\
\hline Rep(\langle(1234)\rangle)& 4 & 15\\
\hline Rep(\langle(12345)\rangle)& 5 & 24\\
\hline Rep(S_3)& 3 & 15\\
\hline Rep(\langle(123)(45)\rangle)& 6 & 15\\
\hline Rep(H_6)& 3 & 12\\
\hline Rep(D_8)& 5 & 12\\
\hline Rep^1(\widetilde{D_8})& 2 & 12\\
\hline Rep(H_{10})& 4 & 12\\
\hline Rep(S_3\times S_2)& 6 & 12\\
\hline Rep^1(\widetilde{S_3\times S_2})& 3 & 12\\
\hline Rep(A_4)& 4 & 14\\
\hline Rep^1(\widetilde{A_4})& 3 & 14\\
\hline Rep(H_{20})& 5 & 9\\
\hline Rep(S_4)& 5 & 8\\
\hline Rep^1(\widetilde{S_4})& 3 & 8\\
\hline Rep(A_5)& 5 & 10\\
\hline Rep^1(\widetilde{A_5})& 4 & 10\\
\hline Rep(S_5)& 7 & 7\\
\hline Rep^1(\widetilde{S_5})& 5 & 7\\
\hline \end{array}
$$

\footnotesize{
{\bf R.B.}: Department of Mathematics,
MIT, Cambridge, MA 02139, USA\\
{\tt bezrukav@math.mit.edu}}

\footnotesize{
{\bf M.F.}: Institute of Information Transmission Problems, and
Independent Moscow University,
Bolshoj Vlasjevskij Pereulok, dom 11,
Moscow 119002 Russia;\\
{\tt fnklberg@gmail.com}}

\footnotesize{
{\bf V.O.}: $\,$  Department of Mathematics,
1222 University of Oregon, Eugene OR 97403-1222, USA\\
{\tt vostrik@math.uoregon.edu}}

\end{document}